\documentclass[a4paper,11pt]{article}
\usepackage{amsmath,amssymb,amsfonts}
\usepackage{graphicx,color}
\usepackage[notcite,notref]{showkeys}
\usepackage{caption}
\usepackage{float}
\usepackage{subcaption}
\usepackage{cite}
\usepackage{verbatim}
\renewcommand{\title}[1]{\begin{center}{\Large \bf #1}\end{center}}
\renewcommand{\author}[2]{\begin{center}{\large #1\\}
{\vspace{0.2cm}\small #2}\end{center}}

\newcommand{\email}[1]{\texttt{#1}}

\newcommand{\ds}{\displaystyle}

\newcommand{\nuu}{s}

\newcommand{\C}{\mathbb{C}}
\newcommand{\HH}{\mathcal{H}}
\newcommand{\bec}{\begin{array}{c}}
\newcommand{\ec}{\end{array}}

\newtheorem{prop}{Proposition}[section]

\newtheorem{rem}{Remark}[section]
\newtheorem{example}{Example}[section]

\date{}

\begin{document}

\title{\Large{\bf{High Order Semi‐implicit Rosenbrock type and Multistep Methods for Evolutionary Partial Differential Equations with Higher Order Derivatives}
}}

\author{{\bf Boscarino Sebastiano, Giuseppe Izzo}}
{
Department of Mathematics and Computer Science, University of Catania (Italy),
\email{sebastiano.boscarino@unict.it},\\
Dipartimento di Matematica e Applicazioni "Renato Caccioppoli", University of Napoli (Italy),
\email{giuseppe.izzo@unina.it}} 


%
%





\begin{abstract}
The aim of this work is to apply a semi-implicit (SI) strategy within a Rosenbrock-type and IMEX linear multistep (LM) framework to a sequence of 1D time-dependent partial differential equations (PDEs) with high order spatial derivatives.  This strategy provides great flexibility to treat {these equations}, and allows the construction of simple linearly implicit schemes without any Newton iteration. Furthermore, the SI schemes so designed 
do not require the severe time-step restrictions typically encountered when using explicit methods
for stability, i.e. $\Delta t = \mathcal{O}(\Delta x^k)$ for the $k$-th order PDEs with $k \ge 2$. For space discretization, this strategy is combined with finite difference schemes. 
We provide examples of methods up to order $p=4$, and
we illustrate the effectiveness of the schemes with  applications to dissipative, dispersive, and biharmonic-type equations. Numerical experiments show that the proposed schemes are stable and achieve 
the expected orders of accuracy.
\\[4mm]

\vspace{3mm}
\end{abstract}

\section{Introduction}
Many time-dependent PDEs which arise in physics or engineering involve the computation
of high order spatial derivatives.
In this paper, we consider a sequence of such PDEs with increasingly higher order derivatives. To simplify the presentation, 
the PDE examples below are restricted to the one-dimensional case.
\begin{itemize}
\item The second order diffusion problem
\begin{equation}\label{First1}
{u_t -} (a(u)u_x)_x = 0,  
\end{equation}
where $a(u) \ge 0$ is smooth and bounded and it is a PDE with second order derivatives. 

Many PDEs of the form (\ref{First1}) which arise in physics or engineering, usually involve the computation of nonlinear diffusion terms, such as: the miscible displacement in porous media \cite{ewing1980galerkin} which is widely used in the exploration of underground water, oil, and gas, the carburizing model \cite{cavaliere2009modeling} which is derived in the chemical heat treatment in mechanical industry, the high-field model in semiconductor device simulations \cite{cercignani2000device, cercignani1997high}, and so on.

In this paper we also consider a one-dimensional version of the convection-diffusion equation
\begin{equation}\label{Conv-Diff}
u_t + f(u)_x - (a(u)u_x)_x = 0. 
\end{equation}

\item The dispersive equation \cite{yan2002local, levy2004local} with third derivatives
\begin{equation}\label{eq:Disper}
u_t + f(u)_x + (r'(u)g(r(u)_x)_x)_x = 0,
\end{equation}
where $f(u)$, $r(u)$ and $g(u)$ are arbitrary (smooth) functions. The Korteweg-de Vries (KdV) equation \cite{korteweg1895xli} which is widely studied in fluid dynamics and plasma physics, is a special case of Eq. (\ref{eq:Disper}) for the choice of the functions $f(u) = u^2$, $g(u) = u$ and $r(u) = u$. The KdV-type equations play an important role in the long-term evolution of initial data \cite{bona1995higher}, are often used to model wave propagation in a variety of nonlinear and dispersive media \cite{benjamin1972model}.

Another choice of the functions $f(u) = u^3$, $r(u) = u^2$  and $g(u) = u/2$, gives the so called general KdV equation \cite{levy2004local}
\begin{equation}\label{eq:Disper2}
u_t + (u^3)_x + (u(u^2_{xx}))_x = 0.
\end{equation}
Eq. (\ref{eq:Disper2}) is known to have compacton solutions of the form:
\begin{equation}
u(x,t) = \biggl\{
\begin{array}{l}
\ds \sqrt{2\lambda} \cos \left(\frac{x-\lambda t}{2}\right), \quad |x-\lambda t|\le \pi,\\[2mm]
0, \quad \quad \textrm{otherwise}.
\end{array}
\end{equation}
Furthermore, Eq. (\ref{eq:Disper2}) is a particular case of the nonlinear dispersive equation \cite{chertock2001particle} 
\begin{equation}\label{eq:Disper2_bis}
u_t + (u^m)_x + (u(u^n_{xx}))_x = 0, \quad m >1, \quad m = n+1,
\end{equation}
with $m = 3$ and $n = 2$. In the numerical tests section we will consider Eq. (\ref{eq:Disper2_bis}) with $m = 3$, $n = 2$ and $m = 2$, $n = 1$. 

Note that in general, the prototype of nonlinear dispersive equations is the $K(m; n)$ equation (\ref{eq:Disper2_bis}), introduced by Rosenau and Hyman in \cite{rosenau1993compactons}. For certain values of $m$ and $n$, the $K(m;n)$ equation has solitary waves which are compactly supported. These structures, the so-called compactons, have several things in common with soliton solutions of the Korteweg–de Vries (KdV) equation where a nonlinear dispersion term replaces the linear dispersion term in the KdV equation.

\item The fourth order diffusion equation
\begin{equation}\label{eq:Fourth}
u_t +  (a(u_x)u_{xx})_{xx} = 0, 
\end{equation}
is a special biharmonic-type equation, where the nonlinearity could be more general but we just present (\ref{eq:Fourth}) as an example. In this paper we will also concentrate on the one dimensional case biharmonic type equation
\begin{equation}\label{eq:Fourth2}
u_t + f(u)_x + (a(u_x)u_{xx})_{xx} = 0.
\end{equation}
 The fourth order diffusion problem has wide applications in the modeling of thin beams and plates, strain gradient elasticity, and phase separation in binary mixtures \cite{han1999dynamics}.
\end{itemize}
For all these equations suitable initial conditions and boundary conditions will be set.

It is well known that the time discretization is a very important issue for time dependent PDEs. For the $k$-th ($k \ge 2$) order PDEs, explicit methods always suffer from stringent and severe time step restriction, i.e., $\Delta t = \mathcal{O}(\Delta x^k)$, for the stability where $\Delta t$ is the time step and $\Delta x$ is the mesh size. This time step is too small, resulting in excessive computational cost and rendering the explicit schemes impractical. 
On the other hand, implicit methods can overcome the drastic time step size restriction imposed by the stability condition for explicit schemes. However, a nonlinear algebraic system must be solved (e.g. by Newton iteration) at each time step. 

In order to apply the SI formulation, 
assume that the semi-discrete formulation of (\ref{Conv-Diff}), (\ref{eq:Disper2_bis}) and (\ref{eq:Fourth2}) can be written as
\begin{equation}\label{DiscS}
\frac{d U(t)}{dt} = \mathcal{F}(U(t)) + \frac{1}{\Delta x^k}\mathcal{B}(U(t))U(t),
\end{equation}

where $U(t)=(U_1(t), U_2(t), ..., U_M(t))^T$, and $U_i(t)$ is the approximate solution at spatial position $x_i$ for $i = 1, ..., M$, i.e.,  $U_i(t) \approx u(x_i,t)$, for $i = 1,...,M$, and uniform grid spacing $\Delta x = x_{i+1}-x_i$, where $\mathcal{F}: \mathbb{R}^M \to \mathbb{R}^M$ and $ \mathcal{B}(U(t))$ is an $M\times M$ matrix. 

As an example, taking the convection nonlinear diffusion equation (\ref{First1}) $\mathcal{F}(U(t))$ is a vector arising from the spatial discretization of $f(u)_x$ and the matrix $\mathcal{B}(U(t)) \in \mathbb{R}^{M \times M}$ is a tridiagonal matrix arising from the discretization of $(a(u)u_x)_x$. Note that the matrix $\mathcal{B}$ inherits its nonlinear 
dependence on $U$ from that of $a(u)$ on $u$.

The semi-implicit (SI) method \cite{boscarino2016high, sebastiano2023high} combines implicit and explicit time discretizations: in the term $\mathcal{B}(U)U$, only the factor $U$ is treated implicitly, while the factor $\mathcal{B}(U)$ and the term $\mathcal{F}(U(t))$ are explicit, thus the resulting scheme can be solved without any iterative solvers. Within the IMEX Runge–Kutta framework, the SI approach is particularly convenient and effective when the problem involves a linearly implicit treatment of the unknown variable in terms with high-order spatial derivatives, as in Eq.~\eqref{DiscS} (see \cite{sebastiano2023high}). Its main advantage is its flexibility, as it can be integrated with various space discretization schemes. 

In this work, we develop an SI strategy based on both a Rosenbrock-type formulation and linear multistep schemes. We present several examples of equations with high-order spatial derivatives of the form \eqref{First1}, \eqref{Conv-Diff}, \eqref{eq:Disper}, \eqref{eq:Fourth}, and \eqref{eq:Fourth2}, which can be efficiently solved using the SI-Rosenbrock-type method and the SI linear multistep approach.

The starting point for the SI Rosenbrock-type formulation is the work in \cite{zhong1996additive}, where the so-called additive semi-implicit Runge–Kutta schemes were introduced and Rosenbrock linearized Runge–Kutta methods were developed. The theoretical framework established in \cite{zhong1996additive} provides a natural basis and can be extended in a straightforward manner to the present semi-implicit setting. 

On the other hand, the starting point of the SI linear multistep approach is the work in \cite{albi2021high}, which develops a general setting for the construction of high order SI linear multistep methods for time-dependent PDEs of the form \eqref{DiscS}, combining an explicit predictor with an implicit corrector. 
However, when the order of the  high-order spatial derivative increases beyond $k \geq 2$ the SI linear multistep scheme tends to become unstable and inaccurate due to the explicit predictor step. To avoid this drawback, we introduce a {\em modified} class of SI linear multistep schemes that turn out to be more efficient, due to a semi-implicit corrector, and more stable and accurate, since it is based on a predictor that is itself semi-implicit (see Sect.~\ref{SIMM}). 
As usual in the context of predictor-corrector (PC) techniques based on linear multistep methods, an order $p^*$ predictor provides a preliminary approximation, and an order $p$ corrector, with $p\geq p^*$, refines it through successive iterations. 
In order to achieve the expected order of convergence $p$ of the scheme, typically $p-p^*$ correction steps are required. 

When applying modified semi-implicit linear multistep (SI-LM) methods, each correction step can improve the accuracy of the solution. In practice, performing approximately $p-p^*$ correction steps per time step, i.e. solving about $p$ linear systems, is sufficient to gradually recover the desired convergence order $p$ of the method.

Finally, we coupled high order finite difference schemes \cite{shu2006essentially} with high-order SI time discretization for solving system (\ref{DiscS}).
We choose the finite difference schemes to discretize the space because of its simplicity in design and coding. For a detailed description of the spatial discretization used in the numerical section, we refer the reader to \cite{sebastiano2023high}.
Nevertheless, other types of spatial discretizations can also be considered, such as local discontinuous Galerkin schemes \cite{cockburn1998local} with application to  various high order PDEs, \cite{levy2004local, shi2019local, tao2020ultraweak, wang2020local, xu2006local, yan2002local}.

 The structure of the paper is as follows. In the next section, we introduce the SI strategy, which is developed using the approach proposed in \cite{boscarino2016high}, partly inspired by partitioned Runge–Kutta methods \cite{wanner1996solving}. Next, in Section \ref{SIRM},  we detail the derivation of the SI-Rosenbrock-type schemes and analyze their stability properties. Several numerical experiments are reported in Section \ref{sect:NR} to evaluate the performance
of the proposed SI-Rosenbrock-type schemes.
In section \eqref{SIMM}, we introduce a {\em modified} class of SI linear multistep (SI-LM) schemes. We analyze basic properties and provide an easy way to construct high-order semi-implicit schemes. Then we provide examples of methods up to order four, based on Backward Differentiation Formulae (BDF) \cite{lambert1991numerical}. 
Finally, Section \ref{NumEx} presents numerical experiments on the test problems considered that confirm the validity of the modified SI-PC-LM approach.

\section{SI methods}\label{SILM}
Starting from system (\ref{DiscS}) we define a multivalued function $\HH:\mathbb{R}^M\times \mathbb{R}^M \to \mathbb{R}^M $ that satisfies $\HH(U,U) = \mathcal{F}(U(t)) + \mathcal{B}(U(t))U(t)$ so that (\ref{DiscS}) can be rewritten as 
\begin{align}
\label{G-systemBIS}
\ds  \frac{dU}{dt}(t) = \HH(U(t),U(t)), \quad \forall\, \, t \geq t_0. 
\end{align}
In \cite{boscarino2016high} the authors introduced a class of SI Runge-Kutta (RK) methods that can be directly applied to (\ref{G-systemBIS}). The SI-RK framework proposed in \cite{boscarino2016high} assumes that function $\HH$ has a non stiff dependence on the first variable and a stiff dependence on the second variable.

Therefore the problem (\ref{G-systemBIS}) can be rewritten as a partitioned form 
\begin{equation}
\left\{
\begin{array}{l}
\ds \frac{d U}{dt}(t) \,=\, \HH(U(t),V(t)),
\\
\,
\\
\ds \frac{dv}{dt}(t) \,=\, \HH(U(t),V(t)),
\end{array}
\right.
\label{ydotzdot}
\end{equation}
with $\HH(U(t),V(t)) =  \mathcal{F}(U(t)) + \mathcal{B}(U(t))V(t)$ and initial conditions $U(t_0) = V(t_0)$,  
where the variable $U(t)$ appearing as the first argument of $\HH$ will be treated explicitly, while $V(t)$ appearing as the second argument will be treated implicitly. In such a case the solution of system (\ref{ydotzdot}) satisfies
  $U(t)=V(t)$ for any $t \ge t_0$ and is also a solution of Equation
(\ref{G-systemBIS}).  System (\ref{ydotzdot}) is a particular case of
partitioned system apparently with an additional computational cost {with respect to (\ref{G-systemBIS})} since we
double the number of variables.

In this section we review the concept of SI RK methods introduced in \cite{boscarino2016high} and next we  derive new classes of SI-Rosenbrock methods.

 A general SI RK method for solving (\ref{ydotzdot}) is a partitioned RK method defined by a double Butcher tableaux with equal weights $b = \tilde{b}$, where the resulting method treats the variable $U(t)$ explicitly using $(\tilde{A}, b)$ and the variable $V(t)$ implicitly using $({A}, b)$. Then the SI RK method applied to (\ref{ydotzdot}) is implemented as follows. Setting $U^{n} = V^n$, we compute the stage values  for $i = 1,...,s$,
\begin{align}\label{IntStag}
U^i = U^n +  \sum_{j = 1}^{i-1} \tilde{a}_{ij}K^j, \quad {V}^i = V^n +\sum_{j = 1}^{i-1} {a}_{ij}K^j + a_{ii} K^i,
\end{align}
and the numerical solution
\begin{align}\label{SolSI}
V^{n+1} = V^n +  {\sum_{i=1}^s} b_i K^i,
\end{align}
where $K_i = \Delta t\HH(U^i,V^i)$, are the RK fluxes.

 We observe that, since $\tilde{b}_i = b_i$, for $i = 1, ...,s$, the numerical solutions coincide, i.e., if $U^0 = V^0$ then $U^{{n}} = V^{n}$,  for all $n\ge 0$. Hence, the duplication of the system is only apparent and not necessary if we adopt the RK fluxes $K_i = \HH(U^{i},V^i)$ as basic unknowns. In this case, the scheme can be equivalently rewritten in the form
 \begin{align}\label{OneS}
K^i = \Delta t\HH(U^{i},  \, U^n + \Delta t \sum_{j = 1}^{i-1} {a}_{ij}K^j + \Delta t a_{ii}K^i),
\end{align}
for $i = 1, ...,s$, with the numerical solution \eqref{SolSI}.
In light of the previous discussion, the main advantage of using the SI RK approach to compute numerical solutions is to solve a linear system. In the case of  system (\ref{DiscS}), Eq. (\ref{OneS}) can be expressed as linear equations in terms of $K_i$:
\begin{equation}\label{Ki}
K^i =\Delta t \mathcal{B}(U^{i})\left( \bar{V}^i+ \Delta t a_{ii}K^i\right),
\end{equation}
where  $ \bar{V}^i = U^n + \sum_{j = 1}^{i-1} {a}_{ij}K^j$ and the matrix {$\mathcal{B}(U^{i})$} are computed explicitly.

\section{SI Rosenbrock-type methods.} \label{SIRM}
A more computationally efficient SI-RK method is a semi-implicit extension of the Rosenbrock-type Runge-Kutta method \cite{wanner1996solving}.  The idea is to consider a semi-implicit scheme (\ref{IntStag})-(\ref{SolSI}) with (\ref{OneS})
and by  linearizing the terms
\begin{equation}\label{KKi}
 K^i =\Delta t \HH(U^i, \bar{V}^i  + \Delta a_{ii}K^i), \; \; \; i=1,\ldots,s,
\end{equation}
with respect to the second variable, we get, 
\begin{equation}\label{Ros}
\left[{\bf I} - \Delta t a_{ii} {\bf J}(U^i, U^n + \Delta t \sum_{j = 1}^{i-1} d_{ij}K^j)\right] K^i \,=\, \, \Delta t\,\HH\left(U^i, \bar{V}_i\right),  
\; \; \; i=1,\ldots,s,
\end{equation}
with $U^i = U^n + \ds\sum_{j = 1}^{i-1}\tilde{a}_{ij}K^j$, and 
\begin{equation*}
\ds V^{n+1} \,=\, V^{n} \,+\,  \sum_{i =1}^{s} {b}_{i}\, K^i, 
\end{equation*} 
where $d_{ij}$ constitute an additional set of parameters,  ${\bf J} = \partial \HH/\partial V \equiv \HH_V$ is the Jacobian matrix with respect to the second argument in (\ref{KKi}), once fixed the first argument. Here {\bf I} denotes the $M \times M$ 
identity matrix. 

  Note that, in the literature, to compute the Jacobian, in classical Rosenbrock methods, it is common to use the same value $a_{ii} = \gamma$ for all $i$ (SDIRK methods). Additionally, $d_{ij} $ are set to 0, so that the Jacobian ${\bf J}( U^i, U^n + \Delta t \sum_{j = 1}^{i-1} d_{ij}k_j)$ is replaced by ${\bf J}(U^n) \equiv {\bf J}(U^n, U^n)$. This approach allows us to calculate the Jacobian only once and use a single LU decomposition to solve equation (\ref{Ros}) for all intermediate stages. 
  
  In particular,  from (\ref{DiscS}), by computing ${\bf J} = \partial \HH / \partial V$, we obtain  ${\bf J} = \mathcal{B}(U(t))$ with ${\bf J}(U^n) = \mathcal{B}(U^n)$. Thus, SI Rosenbrock type methods use the same matrix computed  at the time $t_n$,  namely ${\bf I} - \Delta t \gamma {\bf J}(U^n)$ to solve the linear system for each evaluation of $K^i$, for $i = 1, ...,s$. In contrast, in (\ref{Ki}), SI-RK methods require a recomputing matrix ${\bf I} - \Delta t \gamma \mathcal{B}(U^i)$ to solve the linear system at every stage $K^i$, for every  $i = 1,...,s$.
  
In \cite{zhong1996additive}, several classes of Rosenbrock-type methods are designed for additive semi-implicit Runge-Kutta methods where the function $\HH(U,V)$ has the additive form $\HH(U,V) = f(U) + g(V)$ and this includes cases where $d_{ij} = 0$ and case $d_{ij} = a_{ij}$ with $K^i$ having the additive form as 
$$
K^i = f(U^i) + g(\bar{V}^i + \Delta t a_{ii} K^i),
$$
for $i= 1,...,s$.

Following now the idea in \cite{wanner1996solving}, by introducing on the right hand side an extra term $\Delta t {\bf J}(u^n) \ds\sum_{j = 1}^{i-1} \gamma_{ij}K^j$ and keeping the same Jacobian ${\bf J}(U^n)$ for all stages, with $\gamma_{ii} = \gamma$ and $U^n = V^n$, thus we obtain 
\begin{equation}\label{Ros2}
 K^i \,=\, \, \Delta t\,\HH\left(U^i, \bar{V}_i\right) + \Delta t {\bf J}(U^n) \sum_{j = 1}^{i-1} \gamma_{ij}K^j+  \Delta t \gamma {\bf J}(U^n)K^i,  
\end{equation}
or (\ref{Ros2}) as 
\begin{equation}\label{Ros3}
K^i \,=\, \, \Delta t\left[{\bf I} - \Delta t \gamma {\bf J}(U^n)\right]^{-1} \left(\HH\left(U^i, \bar{V}^i\right) + {\bf J}(U^n)\sum_{j = 1}^{i} \gamma_{ij}K^j\right),  
\end{equation}
with
\begin{equation}\label{ROS_int}
\ds U^i \,=\, U^{n} \,+\,  \sum_{j =1}^{i-1} \tilde{a}_{ij}\, K^j, 
\quad
\bar{V}^i \,=\, V^{n} \,+\, \sum_{j =1}^{i-1} \alpha_{ij} \, K^j, 
\end{equation}
and numerical solution
\begin{equation}\label{RosNS}
\ds V^{n+1} \,=\, V^{n} \,+\,  \sum_{i =1}^{s} {b}_{i}\, K^i.
\end{equation} 
Scheme (\ref{Ros2})-(\ref{ROS_int}) and (\ref{RosNS}) is an $s$-stage SI Rosenbrock method where $\alpha_{ij}$  and $\gamma_{ij}$ are the determining coefficients. Furthermore, for the presentation of order conditions it is convenient to also consider coefficients $\beta_{ij}$ defined by
\begin{equation}\label{Beta}
\beta_{ij} = \alpha_{ij} + \gamma_{ij}, \quad j \le i,
\end{equation}
with $\alpha_{ii} = 0$. A further simplification of the order conditions is possible if we choose $\gamma_{ii} = \gamma$, for all $i$. In the following we will also make use of the abbreviations
 \begin{equation}\label{Beta2}
\tilde{c}_i =\sum_{j = 1}^{i-1} \tilde{a}_{ij}, \quad \alpha_{i} =\sum_{j = 1}^{i-1} \alpha_{ij}, \quad  \beta_{i}' = \sum_{j= 1}^{i-1}\beta_{ij}, \quad \gamma_{i} =\sum_{j = 1}^{i-1} \gamma_{ij} + \gamma, 
\end{equation}
for $i = 1,\dots, s$.

Furthermore, for notational convenience, we introduce the following abbreviations.
The $s$-stage SI Rosenbrock-type scheme is characterized by: an explicit matrix
\begin{equation}\label{BouT}
 \tilde{A} =\left( 
 \begin{array}{ccccc}
   0&0& \dots &0&0\\
  \tilde{c}_2 & 0& \dots&0& 0\\
 \tilde{a}_{31} &\tilde{a}_{32} & \dots& 0&0\\
\vdots & \vdots &  & \ddots & \vdots\\
\tilde{a}_{s1} &\tilde{a}_{s2}&\dots&\tilde{a}_{ss-1}&0
\end{array}\right)
\end{equation}
with $\tilde{b}^T = {b}^T = (b_1, b_2, ..., \gamma)$, $\tilde{c} = (0, \tilde{c}_2, ...,\tilde{c}_s)^T$ and an implicit one
\begin{equation}\label{BouT2}
B = \left(
 \begin{array}{ccccc}
   \gamma &0& \dots &0&0\\
  \beta_{21} & \gamma & \dots&0& 0\\
 \beta_{31} &\beta_{32} & \gamma& 0&0\\
\vdots & \vdots &  & \ddots & \vdots\\
\beta_{s1} &\beta_{s2}&\dots& \dots& \gamma
\end{array}\right).
\end{equation}
\begin{rem}
Note that in the case of non-autonomous systems, i.e., $U' = \HH(t,U,U)$ in the computation of the Jacobian and an extra term appears, i.e. $\partial{\HH}/\partial{t}$, and then the order conditions change and become significantly more complicated for $p\ge 3$. From now on, we restrict our analysis of the order conditions to autonomous systems up to third-order. For further details on the non-autonomous case, we refer the reader to \cite{wanner1996solving}.
\end{rem}
 
{\bf Order conditions up to third-order SI-Rosenbrock-type method.} Here, following the theory developed in \cite{wanner1996solving}, we derive the order conditions on the free parameters, which ensure that the method is of a given order $p$. Specifically, considering the numerical solution after one time step,
$$
V^1 = V^0 + \sum_{i = 1}^sb_i K^i,
$$
the local error satisfies
$$
V(t_0 + \Delta t) - V^1 = \mathcal{O}(\Delta t^{p+1}).
$$
As a first step, we compute the  $q$-derivatives of the numerical solution at $h =0$:
\begin{equation}\label{OrdCond}
V^{1,(q)} = \sum_{i=1} b_i K_i^{(q)}(0).
\end{equation}
Hence, it remains to compute the derivatives of $K_i$ from (\ref{OrdCond}) with respect to $\Delta t$ evaluated at $\Delta t = 0$. By applying Leibniz’s rule, this leads to
\begin{equation}\label{KK_i}
K_i^{(q)}|_{\Delta t= 0} = q(\HH(U^i,\bar{V}^i))^{(q-1)}|_{\Delta t= 0} + q {\HH}_{V}^0\sum_{j = 1}^{i}\gamma_{ij} k_j^{(q-1)}|_{\Delta t= 0}.
\end{equation}
with ${\HH}_{V}^0 \equiv {\bf J}(U_0) $.
Now we compute the order conditions up to third order. 
We start by considering the system
\[
U' = \HH(U,V), \qquad V' = \HH(U,V),
\]
and we compute the following derivatives up to third order. For the second and third derivatives, we obtain
\[
U'' = V'' = \HH_{U} \, U' + \HH_{V} \, V' 
       = \HH_{U} \, \HH + \HH_{V} \, \HH,
\]
and
\begin{equation}\label{Ex3}
\begin{aligned}
U''' = V''' = \; & \HH_{UU}(\HH,\HH) + \HH_{VV}(\HH,\HH) 
+ 2\,\HH_{UV}(\HH,\HH) \\
& + \HH_{U}\HH_{U}\HH + \HH_{U}\HH_{V}\HH 
+ \HH_{V}\HH_{U}\HH + \HH_{V}\HH_{V}\HH.
\end{aligned}
\end{equation}
where \(\HH_{U}\),\(\HH_{V}\), \(\HH_{UU}\), \(\HH_{VV}\) and \(\HH_{VU} = \HH_{UV}\) denote the partial derivatives of \(\HH\) with respect to \(U\) and \(V\), respectively. 

Now, from (\ref{KK_i}) we compute:
$$
K_i^{(0)}|_{\Delta t= 0} = 0, \quad K_i^{(1)}|_{\Delta t= 0} = \HH(U_0,V_0) \equiv  \HH^0,
$$
$$
K_i^{(2)}|_{\Delta t= 0}  = 2 \sum_{j= 1}^{i-1}\tilde{a}_{ij} \HH_{U}^0 \HH^0 + 2 \left( \sum_{j = 1}^{i-1}\alpha_{ij} +  \sum_{j = 1}^{i} \gamma_{ij} \right){\HH}_{V}^0 \HH^0,
$$
or, equivalently,
$$
K_i^{(2)}|_{\Delta t= 0}  = 2 \sum_{j= 1}^{i-1}\tilde{a}_{ij} \HH_U^0 \HH^0 + 2 \left( \sum_{j = 1}^{i-1}\beta_{ij} + \gamma \right){\HH}^0_{V} \HH^0.
$$
Here, the subscript index $0$ indicates that the derivatives are evaluated at $(U_0, V_0)$. 

Then, inserting into (\ref{OrdCond}) and comparing the derivatives of the numerical solution 
(for $q \geq 1$) with those of the exact solution $V'$ and $V''$, we obtain the following 
order conditions for the second order:
\begin{equation}\label{Cond_2}
\sum_{i = 1}^s b_i = 1, \quad  \sum_{i= 1}^{s}b_i\tilde{c}_{i} = \frac{1}{2}, \quad \sum_{i = 1}^s b_i\beta_{i}' = \frac{1}{2} - \gamma.
\end{equation}
Here we used the abbreviations
$c_i = \sum_{j=1}^{i-1} \tilde{a}_{ij}$, for $i = 1,...,s$.

Now we derive conditions up to third order. From (\ref{KK_i}), and setting $q = 3$, we obtain
$$
K_i^{(3)}|_{\Delta t= 0} = 3(H(U^i,\bar{V}^i))^{(2)}|_{\Delta t= 0} + 3 {\HH}^0_{V}\sum_{j = 1}^{i}\gamma_{ij} K_j^{(2)}|_{\Delta t= 0}.
$$
where
$$
3(H(U^i,\bar{V}^i))^{(2)}|_{\Delta t= 0} = 3 \left(H_{UU}({U}',{U}') + H_{VV}(\bar{V}',\bar{V}') + 2H_{UV}({U}',\bar{V}') + H_U{U}'' + H_V{\bar{V}}'' \right)|_{\Delta t= 0},
$$
with $U^{(q)}|_{\Delta t = 0} = \sum_{j = 1}^{i-1} \tilde{a}_{ij}K_j^{(q)}|_{\Delta t = 0}$, $\bar{V}^{(q)}|_{\Delta t = 0} = \sum_{j = 1}^{i-1} \alpha_{ij}K_j^{(q)}|_{\Delta t = 0}$.

Furthermore,
$$
 3\HH^0_{V} \sum_{j = 1}^{i}\gamma_{ij} K_j^{(2)}|_{\Delta t= 0} = 6\left(\sum_{j = 1}^{i}\gamma_{ij}\tilde{a}_{ij}\HH_{V}^0\HH_{U}^0 \HH^0 + \sum_{j = 1}^{i}\gamma_{ij}\left( \sum_{k = 1}^{j-1}\beta_{jk} + \gamma\right)\HH_{V}^0 \HH_{V}^0 \HH^0 \right).
 $$
 Comparing the derivatives of the numerical solution (\ref{OrdCond}) with $q = 3$, and  those of the exact solution (\ref{Ex3})
we have the following third-order conditions:
 \begin{equation}\label{Cond3_1}
\sum_{i,j}b_i\tilde{c}^2_{i} = \frac{1}{3}\quad \sum_{i,j}b_i\tilde{a}_{ij}\tilde{c}_{j} = \frac{1}{6},
 \end{equation}
\begin{equation}\label{Cond3_2}
\sum_{i,j,k}b_i\tilde{a}_{ij}\alpha_{ik} = \frac{1}{3}, \quad \sum_{i,j,k} b_i\alpha_{ij}\alpha_{ik} = \frac{1}{3},
\end{equation}
\begin{equation}\label{thirdB}
\sum_{i,j,k}b_i\tilde{a}_{ij}(\alpha_{jk} + \gamma_{jk} )= \frac{1}{6}, \quad
\sum_{i,j,k}b_i(\alpha_{ij} + \gamma_{ij} )\tilde{a}_{jk}= \frac{1}{6},\quad
\sum_{i,j,k}b_i(\alpha_{ij} + \gamma_{ij} )(\alpha_{jk} + \gamma_{jk} )= \frac{1}{6}.
\end{equation}
Now using (\ref{Beta2}) the formulas (\ref{Cond3_2}) and (\ref{thirdB}) become
\begin{align}\label{thirdBbis}
\begin{split}
\sum_{i,j}b_i\tilde{a}_{ij}\alpha_{i} &= \frac{1}{3}, \quad \sum_{i,j}b_i \alpha_i^2 = \frac{1}{3},\\
\sum_{i,j}b_i\tilde{a}_{ij}\beta'_{j}&= \frac{1}{6}- \frac{1}{2}\gamma, \quad
\sum_{i,j,k}b_i\beta_{ij}\tilde{c}_{j}= \frac{1}{6}- \frac{1}{2}\gamma,\quad
\sum_{i,j}b_i\beta_{ij} \beta'_{j}= \frac{1}{6} - \gamma + \gamma^2.
\end{split}
\end{align}

{\bf Stability function.}
We consider the simplified one-dimensional linear model equation of the form (\ref{G-systemBIS}) 
\[
U' = \lambda U + \mu U,
\] 
with $\lambda, \, \mu \in \mathbb{C}$, and the  corresponding partitioned form (\ref{ydotzdot}) becomes:
\begin{equation}\label{PAs}
U' = \lambda U + \mu V, \quad V' = \lambda U + \mu V,
\end{equation}
with 
\[
\mathcal{H}(U,V) = \lambda U + \mu V, \quad U(0) = V(0) = 1.
\]

Here $\HH_V^0 = {\bf J}(U^0) = \mu$, and applying the method (\ref{Ros2})-(\ref{ROS_int}) and (\ref{RosNS}) to \eqref{PAs}, the numerical solution after one step reads $U^1 = R(\tilde{z}, z)$ with 
\begin{equation}\label{Rfun}
 R(\tilde{z}, z) = 1 + (\tilde{z} + z) b^T(I - \tilde{z}\tilde{A} - z B)^{-1}{\bf e},
\end{equation}
where $b^T = (b_1, ..., b_s)$ and $B = (\beta_{ij})_{i,j=1}^s$. The stability region of the method is defined as
\begin{equation}\label{Sreg}
S = \{ (\tilde{z},z): \tilde{z}, z \in \mathbb{C} \quad \textrm{such that} \quad |R(\tilde{z},z)|<1\} \subset \mathbb{C}^2.
\end{equation}
Note that, in \eqref{Rfun}, for $z \to \infty$
\begin{equation}\label{Rz}
R(\tilde{z}, \infty) = 1 -  b^T B^{-1}{\bf e} = 1 - \sum_{i,j} b_i \omega_{ij},
\end{equation}
where we denote with  $\omega_{ij}$ the coefficients of the inverse of the matrix $(\beta_{ij})$. 

{Whenever the coefficients of the Rosenbrock method are chosen such
that \cite{wanner1996solving}, 
\begin{equation}\label{SA}
\beta_{si} = \alpha_{si} + \gamma_{si} = b_i, \quad i = 1, ..., s, \quad \textrm{and} \quad \alpha_s =\sum_{j=1}^{s-1}\alpha_{sj} = 1, 
\end{equation}
they are called {\em stiffly accurate}. Since in this case $b^T = e^T_sB$, (or $b^TB^{-1} =e^T_s$) with  $e^T_s = (0,...,0,1)$, from \eqref{Rz} we get $R(\tilde{z},\infty) = 0$, that is, for 
fixed $\tilde{z} \in \C$, this ensures $L$-stability for the Rosenbrock method.
}
Note that conditions (\ref{SA}) together with $\sum_{i}^s b_i = 1$, implies $\gamma_s = \sum_{j}^{s-1} \gamma_{sj} + \gamma = 0$.

{The SI-Rosenbrock type methods with $\gamma>0$  possess the stability functions  of the form
\begin{equation}\label{SFR}
R(\tilde{z},z) = \frac{P(\tilde{z}, z)}{(1-\gamma z)^s},
\end{equation}
with $P(\tilde{z}, z)$ polynomial of degree $\ge s$ similar to the IMEX-RK methods \cite{boscarino2024implicit}.
For instance, the one-stage method, \eqref{Ros2}-\eqref{RosNS} reads as:
$$ V^{n+1} = V^n + K^1, \quad K^1 = (\tilde{z} + z)V^n + \gamma z K^1,
$$
with stability function
$$
R(\tilde{z},z) = \frac{1 + \tilde{z} + (1 - \gamma)z}{1-\gamma z}.
$$ 
For $\gamma = 1$ we obtain the stability function of the IMEX Euler method. Similarly, for the two-stage SI Rosenbrock-type method \eqref{Ros2}--\eqref{RosNS}, 
satisfying the order conditions \eqref{Cond_2}, namely
\[
b_1 = 1-b_2, \quad b_2 = \tfrac{1}{2\tilde{c}_2}, \quad b_2 \beta_{21}' = \tfrac{1}{2} - \gamma,
\]
with $b_2 \neq 0$ and for any choice of $\gamma > 0$, the method is of order two.  
The corresponding stability function is the rational function \cite{boscarino2024implicit}
\[
R(\tilde{z}, z) = 
\frac{\tfrac{1}{2}\tilde{z}^2 + \tilde{z} + 1 \;+\; z(1-2\gamma)(1+\tilde{z}) \;+\; z^2\left(\tfrac{1}{2} - 2\gamma + \gamma^2\right)}
{(1 - z\gamma)^2},
\]
which coincides with that of the second-order two-stage IMEX-RK scheme.  
By this analogy, we choose in the following $\gamma \geq \tfrac{1}{4}$ in order to guarantee $A$-stability for the method, \cite{boscarino2024implicit}.}

Finally, to understand the geometrical structure of the region $S$ \eqref{Sreg}  associated with a SI Rosenbrock-type method, we follow the idea in \cite{pareschi2000implicit, boscarino2024implicit}. Then  we define the region $\tilde{S}$ as the largest stability region for $\tilde{z} \in \C$ such that the scheme is still $A$-stable. Such region is formally defined as
$$
\tilde{S} = \{ \tilde{z} \in \mathbb{C}:  \quad |R(\tilde{z},z)|\le 1, z \in \mathbb{C}^{-}\} .
$$
Explicit expressions for the set $\tilde{S}$ have been given in \cite{pareschi2000implicit}.
In particular, to compute $\tilde{S}$, we consider $z = iy$, $y \in \mathbb{R}$, and use the maximum modulus principle, then we can rewrite $S$ as
$$
\tilde{S} = \{ \tilde{z} \in \mathbb{C} \quad \textrm{such that} \quad \max_{y \in \mathbb{R}}|R(\tilde{z},iy)|\le1\} ,
$$
and the boundary of this region (boundary locus) is obtained as
\begin{equation}\label{Stab_region}
\partial \tilde{S} = \{ \tilde{z} \in \mathbb{C} \quad \textrm{such that} \quad \max_{y \in \mathbb{R}}|R(\tilde{z},iy)|=1\} .
\end{equation}
Next we show stability regions $\tilde{S}$ of third-order SI-Rosenbrock-type methods that we will introduce in the following section.

{\bf Construction of methods of order 3.}
In order to construct $s$-stage SI Rosenbrock-type methods of order 3, we need to require $s = 4$. The case $s = 3$ results in a system of equations that is unsolvable with respect to the number of coefficients that need to be determined for our method.  

We list, for convenience, the whole set of order conditions (\ref{Cond_2})-(\ref{Cond3_1})-(\ref{Cond3_2})-(\ref{thirdBbis}) for the case $s = 4$:
\begin{align}\label{2}
\begin{split}
&b_1  = 1 - b_2 - b_3  - \gamma,\\
&b_2 \tilde{c}_2 + b_3 \tilde{c}_3 + \gamma \tilde{c}_4 = 1/2,\\
&b_2 \beta'_2 + b_3 \beta'_3 + \gamma \beta'_4 = 1/2 - \gamma,\\  
\end{split}
\end{align}
\begin{align}\label{3}
\begin{split}
&b_2 \tilde{c}_2^2 + b_3 \tilde{c}_3^2 +  \gamma \tilde{c}_4^2 = 1/3,\\
&b_3 \tilde{a}_{32}\tilde{c}_2 +  \gamma \left(\tilde{a}_{42}\tilde{c}_2 + \tilde{a}_{43}\tilde{c}_3\right)= 1/6,\\
&b_3 \tilde{a}_{32}\alpha_2 + \gamma \left(\tilde{a}_{42}\alpha_2+\tilde{a}_{43}\alpha_3\right)  = 1/3,\\
&b_2 \alpha_2^2 + b_3 \alpha_3^2   = 1/3 -\gamma,\\
&b_3\tilde{a}_{32}\beta'_2 + \gamma\left(\tilde{a}_{42}\beta'_2 + \tilde{a}_{43}\beta'_3\right) =1/6 - \gamma/2,\\
&b_3 \beta_{32} \tilde{c}_2 + \gamma\left( \beta_{42} \tilde{c}_2 + \beta_{43} \tilde{c}_3\right)= 1/6 - \gamma/2,\\
& b_3\beta_{32}\beta'_2 + \gamma\left(\beta_{42}\beta'_2 + \beta_{43}\beta'_3\right)  = 1/6 - \gamma + \gamma^2.
\end{split}
\end{align}
We impose conditions \eqref{SA}:
\begin{align}\label{4}
b_i = \beta_{4i}, \quad \alpha_4 = 1, \quad i = 1, ..., 4,
\end{align}
such that the method is SA, i.e., $R(\infty) = 0$, with $b_4 = \gamma$. Note that from (\ref{4}), the sixth equation in (\ref{3}) can be rewritten as
$$
b_3 \beta_{32} \tilde{c}_2 + \gamma\left( b_2\tilde{c}_2 + b_3\tilde{c}_3\right)= 1/6 - \gamma/2,
$$
and by the the second equation in (\ref{2}) we get 
\begin{align}\label{5}
b_3 \beta_{32} \tilde{c}_2 + \gamma\left( 1/2 - \gamma \tilde{c_4}\right)= 1/6 - \gamma/2.
\end{align}
Since the number of unknown coefficients to be determined is larger than the number of equations, the system has several degrees of freedom. To reduce that, we impose the additional condition
\[
\sum_i b_i \alpha_i = \tfrac{1}{2},
\] 
which ensures that the pairs $(b_i, \alpha_i)$ corresponds to the coefficients of a standard quadrature formula \cite{wanner1996solving}. As a consequence, we include among the conditions above the extra relation
\begin{align}\label{6}
b_2 \alpha_{2} + b_3 \alpha_{3} = \tfrac{1}{2} - \gamma.
\end{align}

{The computation of the coefficients $\alpha_{ij}$, $\beta_{ij}$, $\gamma$ and $b_i$ is now done by considering the order conditions (\ref{2}), (\ref{3}), (\ref{5}), and (\ref{6}). 
\\
The main goal of this paper is to design a third-order SI Rosenbrock scheme with good stability properties, in particular $L$-stability. The procedure outlined below is not the only possible way to determine the coefficients.  Several free parameters can be chosen in different ways to further improve stability or to minimize the error constant. A detailed optimization of these aspects, however, lies beyond the scope of the present paper.
\\
{\em Step 1.} We choose $\gamma >0$ such that the stability function (\ref{Rfun}) has the desired stability properties. In particular, as mentioned above, we take $\gamma \ge 1/4$. At this stage, we regard $b_2$, $\alpha_2$, and $\tilde{a}_{32}$ as free parameters.}\\
{{\em Step 2.} By setting $b_2 = 0$, and solving the fourth equation in (\ref{3}) together with equation (\ref{6}) we obtain $b_3$, and $\alpha_3$ as well as $b_1 = 1- b_2 - b_3 - \gamma$.\\
{\em Step 3.} From the second equation in (\ref{2}) and the second equation in (\ref{3}), we solve for $\tilde{c}_2$ and $\tilde{c}_3$ in terms of $\tilde{c}_4$, $\tilde{a}_{42}$ and $\tilde{a}_{43}$. Substituting these expressions into the first equation of (\ref{3}), we obtain two possible values for $\tilde{c}_4$. By selecting one of them $\tilde{c}_4$, we then determine  $\tilde{c}_2$ and $\tilde{c}_3$, which depend only on $\tilde{a}_{42}$ and $\tilde{a}_{43}$.\\
{\em Step 4.} From the third equation in \eqref{2} and  using $\beta'_4 = b_1 + b_2 + b_3 = 1-\gamma$  we compute $\beta_3'$.
Substituting this value into the seventh equation of \eqref{3}, we obtain $\beta_2'$ which depends on $\beta_{32}$.\\
{\em Step 5.} Setting $\tilde{a}_{32} = 0$ and  considering that $\tilde{c}_2$ depends on $\tilde{a}_{42}$ and $\tilde{a}_{43}$, we use the  third equation in \eqref{3} together with  equation \eqref{5} to compute $\tilde{a}_{42}$ and $\tilde{a}_{43}$ in terms of $\beta_{32}$.\\ 
{\em Step 6.} Substituting the expressions for $\tilde{a}_{42}$, $\tilde{a}_{43}$ and $\beta_2'$ and $\beta_3'$ into the fifth equation of \eqref{3}, we determine $\beta_{32}$.
\\
{\em Step 7.} Once $\beta_{32}$ is known, we can explicitly compute $\beta_2'$, $\tilde{a}_{42}$, $\tilde{a}_{43}$, Consequently, this also fixes the values of $\tilde{c}_2$ and $\tilde{c}_3$.\\
{\em Step 8.} The remaining coefficients can then be obtained from formulas \eqref{Beta}, with  \eqref{Beta2} and \eqref{4}.\\ 
}
Next, we consider different SI Rosenbrock-type schemes by selecting specific values of $\gamma$. In particular we consider the values: 
$$\gamma =1-1/\sqrt{2}, \quad \gamma =13/50, \quad \gamma = 3/4.
$$

As an example, 
we present the SI Rosenbrock scheme with $\gamma = 3/4$.

From step 1. we have $\gamma = 3/4$ and $b_2 = \tilde{a}_{32} = 0$ with $\alpha_2 = \alpha_{21} = 2 \gamma = 3/2$. 

From step 2. we obtain
$$ b_3 = -3/20, \quad \alpha_3 = 5/3,
$$
with $b_1 = 2/5$. Since $\alpha_{31} + \alpha_{32} = \alpha_3$ by setting $\alpha_{31} = 0$ we get $\alpha_{32} = \alpha_3$. 

From step 3. we choose $\tilde{c}_4 = 1$, which gives $\tilde{c}_3 = 5/3$. With $\tilde{a}_{32} = 0$ we have $\tilde{a}_{31} = \tilde{c}_3$. Using $\beta_4' = 1-\gamma$, step 4. yields $\beta_3' = 35/12$. 

From steps 5. and 6.  we obtain $\beta_{32} = 65/108$, and from step 7., we get 
$$
\tilde{a}_{42} = 53/297, \quad \tilde{a}_{43} = 6/55, \quad \tilde{c}_2 = 3/13 \quad \beta_2' =\beta_{21} = -177/52.
$$ 

Finally, step 8. gives the remaining coefficients. First, 
$$
\gamma_{21} = \beta_{21} - \alpha_{21} = \beta_2' - \alpha_2 = -255/52,$$ 
with 
$$\tilde{a}_{41} = \tilde{c}_4 - \tilde{a}_{42} - \tilde{a}_{43} = 1063/1485,
$$ and 
$$
\beta_{31} = \beta_3'-\beta_{32} = 125/54.
$$
Thus,
$$\gamma_{31} = \beta_{31} - \alpha_{31} =\beta_{31} = \beta_3' - \beta_{32} = 125/54,\quad 
 \gamma_{32} = \beta_{32} - \alpha_{32} = -115/108.$$

Finally, since 
$$
\alpha_4 = \sum_{i = 1}^3 \alpha_{4i} = \sum_{i = 1}^3 (\beta_{4i} -\gamma_{4i}),
$$ 
by setting $\alpha_{41} = \alpha_{43} = 0$ we obtain
$$ 
\gamma_{41} = b_1 = 2/5, \quad \gamma_{43} = b_3 = -3/20.
$$
Moreover, since
$$
\alpha_{42} = \beta_{42} -\gamma_{42} = b_2 -\gamma_{42} = -\gamma_{42}= \tilde{c}_4 = 1,
$$
we deduce $\gamma_{42} = -1$.
Below, we explicitly present the scheme together with the computed coefficients:
\begin{align}\label{SIR4}
\begin{split}
{\bf A} &= {\bf I} - \frac{3}{4} \Delta t {\bf J}(U^n),\\
{\bf A} K^1 &= \Delta t \HH(U^n, V^n),\\
{\bf A} K^2 &= \Delta t \left(\HH(U^n + \frac{3}{13}K^1, V^n + \frac{3}{2}K^1) -{\bf J}(U^n)\frac{ 255}{52}K^1\right),\\
{\bf A} K^3 &= \Delta t \left(\HH(U^n + \frac{5}{3}K^1 , V^n + \frac{5}{3}K^2) + {\bf J}(U^n) \frac{125}{54}K^1 - {\bf J}(U^n)\frac{115}{108} K^2\right),\\
{\bf A} K^4 &= \Delta t \left(\HH(U^n + \frac{1063}{1485}K^1 + \frac{52}{297}K^2 + \frac{6}{55}K^3, V^n + K^2) + {\bf J}(U^n)(\frac{2}{5}K^1 - K^2 - \frac{3}{20}K^3) \right),\\
V^{n+1} &= V^n + \frac{2}{5}K^1 - \frac{3}{20}K^3 +  \frac{3}{4} K^4. 
\end{split}
\end{align}

\begin{figure}[H]
\begin{center}
\includegraphics[scale=0.4]{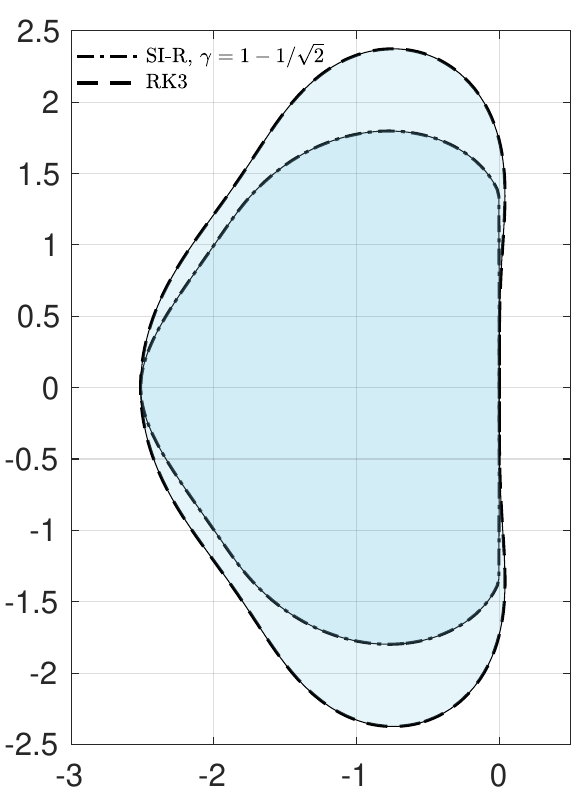}\quad
\includegraphics[scale=0.4]{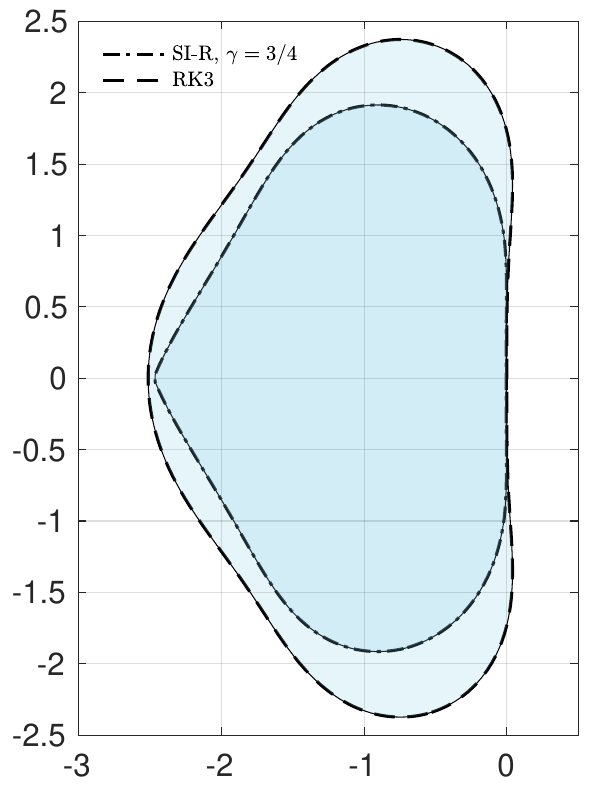}\quad
\includegraphics[scale=0.4]{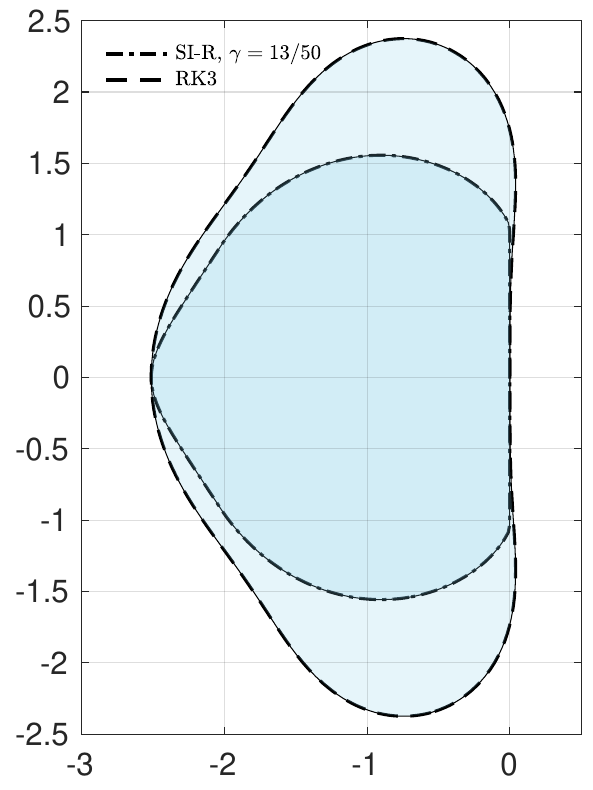}
\end{center}
\caption{Stability region (\ref{Stab_region}) for different values $\gamma$. Left $\gamma = 1-1/\sqrt{2}$, middle $\gamma = 3/4$, right $\gamma = 13/50$.}
\label{fig_stab}
\end{figure}

Finally, we show in Figure \ref{fig_stab} the stability regions $\tilde{S}$ of these SI Rosenbrock type schemes with different values of $\gamma$.

\subsection{Numerical  Results.} \label{sect:NR}
Below, we will numerically verify the orders of accuracy and performance of the third order SI-Rosenbrock type schemes for different values of $\gamma$ for the solution of time dependent PDEs with high order derivatives: (1), (2), (6), (7) and (8). For simplicity we used in all of our examples periodic boundary conditions, although most of our discussions can be adapted for other types of boundary conditions.

For the space discretization we use high-order finite difference schemes introduced in \cite{shu2006essentially, sebastiano2023high}, and we discretize the  convection term  by means of standard third order finite difference WENO schemes with local Lax-Friedrichs flux,  \cite{jiang1996efficient, liu1994weighted, shu2006essentially}. All the tests presented in this section are performed in one-dimension since the main issue of this paper is to design an appropriate high order time discretization for time-dependent PDEs.

\begin{itemize}
\item[Test 1.] Next we consider the convection-diffusion equation
\begin{equation}\label{Eq:2Diff}
u_t + (u^2/2)_x = (a(u)u_{x})_{x} + f(x,t), \quad x\in (-\pi, \pi)
\end{equation}
where the diffusion coefficient is $a(u) = u^2 + 2$, initial condition $u(x,0) = \sin(x)$ and the source term
$$
f(x,t) = \frac{1}{4}\left( 4\cos(x+t) + 9\sin(x+t) + 2\sin(2(x+t)) - 3 \sin(3(x+t))\right).
$$
The problem has the exact solution: $u(x,t) = \sin(x+t)$. 

The standard third order WENO(3,2)  scheme \cite{shu2006essentially}  
is used for the discretization of the convection term $(u^2/2)_x$. We compute to $T = 1$ with the time step $\Delta t =  \Delta x$. The numerical results are listed in Table \ref{table:2}. 
These results confirm the stability and the expected order of accuracy of the SI-Rosenbrock methods coupled with finite difference schemes.

	\begin{table}[htp]
\begin{center}
\begin{tabular}{c|c|c|c|c|c|c|c}
\hline
Scheme     & N     & $L^2$-error & order& $L^1$-error & order& $L^\infty$-error & order \\
\hline
SI-R, $\gamma = 13/50$ &  40   &  5.7907e-04 & - &      5.4529e-04   &   - & 9.1249e-04 &     -     \\
                                       &  80   &     8.8345e-05   & 2.71 &8.8124e-05& 2.63& 1.1886e-04 & 2.94\\
                                       &  160   &    1.2787e-05 & 2.78& 1.2933e-05& 2.76& 1.4922e-05&3.00\\
                                      &  320   &       1.8538e-06  & 2.79 & 1.8808e-06 & 2.78 &2.0453e-06&2.87\\
                                      &  640   &       2.4342e-07  & 2.93  & 2.4629e-07& 2.93& 2.6306e-07& 2.96 \\
\hline
SI-R, $\gamma = 3/4$   &  40   &  3.8819e-03 & - &      3.6931e-03   &   - & 5.2543e-03 &     -     \\
                                       &  80   &   7.9215e-04  & 2.30 &7.6735e-04& 2.26& 1.0492e-03 & 2.32\\
                                       &  160   &  1.4133e-04 & 2.49& 1.3994e-04& 2.46& 1.7398e-04&2.60\\
                                      &  320   &       2.4209e-05  & 2.55 & 2.4271e-05 & 2.53 &2.8472e-05&2.61\\
                                      &  640   &     3.7624e-06  & 2.68  & 3.7299e-06& 2.70& 4.4558e-06& 2.67 \\
                                      \hline
SI-R, $\gamma = 1-1/\sqrt{2}$   &  40   & 6.1483e-04 & - &         5.7829e-04   &   - & {9.7074}e-04 &      -    \\
                                                    &  80   &   9.6939e-05  & 2.66 &9.7006e-05& 2.58& 1.2706e-04 &2.93\\
                                                   &  160   &   1.4473e-05 &2.74& 1.4714e-05&2.72& 1.6298e-05&2.96\\
                                                   &  320   &        2.163e-06  & 2.74 & 2.1909e-06 &2.75 &2.4405e-06&2.74\\
                                                   &  640   &   2.8932e-07  & 2.90  & 2.9151e-07& 2.90& 3.2393e-07& 2.91\\
                                                   \hline
\end{tabular}
\end{center}
\caption{Test 1. The $L^2$, $L^1$, $L^{\infty}$ errors and orders of accuracy for equation (\ref{Eq:2Diff}).}
\label{table:2}
\end{table}%
\item[{Test 2.}] 
In this test, we check the accuracy of the SI-Rosenbrock type scheme applied to Eq. (\ref{eq:Disper2})  {\cite{tan2022stability}}
\begin{equation}\label{KdV_22}
u_t +  (u^3)_x +(u(u^2)_{xx})_x = 0, \quad x \in \left(-\frac{3}{2}\pi, \frac{5}{2}\pi\right),
\end{equation}
with initial condition $u(x,0) = \sqrt{2 {\lambda}}\cos(x/2)$ and exact solution 
$$
u(x,t) = \sqrt{2 {\lambda}}\cos((x-{\lambda} t)/2).
$$

We compute to $T = \pi$, with ${\lambda} = 0.1$. We choose $\Delta t = \Delta x$.  In Table \ref{table:KdV} we show the convergence rate for these schemes. As expected, the SI Rosenbrock type schemes are stable and achieve the correct order of accuracy. 

Note that, since $\lambda$ represents the velocity of the traveling wave, selecting a larger value of $\lambda$ requires a careful choice of the time step to ensure the stability of the method. In such cases, we set the time step as 
\begin{equation}\label{Delta_t} 
\Delta t = CFL \frac{\Delta x}{\max_{u} |f'(u)|}
\end{equation}
with a fixed $CFL$.

As an example, we provide Table \ref{table:KdV_10} below that  shows the  convergence rates of the schemes SI Rosenbrock type for different values of $\gamma$  with $\lambda = 10$,  $CFL = 0.5$ and $T = \pi/4$.  We see that the schemes are stable and achieve the expected order of accuracy. Note that larger $\lambda$ causes larger errors.  
	\begin{table}[htp]
\begin{center}
\begin{tabular}{c|c|c|c|c|c|c|c}
\hline
Scheme     & N     & $L^2$-error & order & $L^1$-error & order & $L^{\infty}$-error & order  \\
\hline
SI-R, $\gamma = 3/10$ &  80   &  8.3717e-04 & - &      7.4284e-04 &   - &1.3899e-03& -         \\
                                       &  160   &     1.0634e-04 & 2.97& 9.289e-05& 3.00&2.6686e-04&2.38\\
                                      &  320   &       1.3408e-05 & 2.98 & 1.1537e-05 & 3.00 &4.4896e-05&2.57\\
                                      &  640   &       1.6915e-06  & 2.98  & 1.4397e-06& 3.00& 7.3241e-06& 2.62  \\
\hline
SI-R, $\gamma = 3/4$   
                                       &  80   &   8.4776e-04  &- & 7.4989e-04& - & 1.4355e-03& - \\
                                       &  160   &  1.0768e-04 & 2.97& 9.3767e-05& 3.00& {2.7374}e-04&2.39\\
                                      &  320   &      1.3578e-05  & 2.98 &  1.1653e-05 & 3.00&4.5897e-05&2.58\\
                                      &  640   &   1.7128e-06  & 2.98  &1.4538e-06& 3.00&  7.4797e-06& 2.62 \\
                                      \hline
SI-R, $\gamma = 1-1/\sqrt{2}$  
                                                    &  80   &   8.3718e-04  & -& 7.4284e-04& -& 1.3900e-03 &-\\
                                                   &  160   &    1.0634e-04 &2.98& 9.289e-05&3.00& 2.6686e-04&2.38\\
                                                   &  320   &      1.3408e-05  &2.99 & 1.1537e-05 & 3.00&4.4896e-05& 2.57 \\
                                                   &  640   &  1.6915e-06 & 2.99  &  1.4397e-06& 3.00& 7.3242e-06&2.62\\
                                                   \hline
\end{tabular}
\end{center}
\caption{Test 2. The $L^2$, $L^1$, $L^{\infty}$ errors and orders of accuracy for equation (\ref{KdV_22}).}
\label{table:KdV}
\end{table}
\begin{table}[htp]
\begin{center}
\begin{tabular}{c|c|c|c|c|c|c|c}
\hline
Scheme     & N     & $L^2$-error & order & $L^1$-error & order & $L^{\infty}$-error & order  \\
\hline
SI-R, $\gamma = 3/10$ &  80   &  2.9585e-03 & - &      2.7744e-03 &  - &4.4434e-03&  -        \\
    &  160   &   4.0981e-04& 2.85& 3.240e-04& 3.09 &  7.5084e-04&2.56\\            &  320   &        5.7821e-05 & 2.82 & 4.0098e-05 & 3.01 &1.7626e-04&2.09\\               &  640   &7.3668e-06  & 2.97  & 4.9762e-06& 3.01& 2.3916e-05& 2.88  \\
\hline
SI-R, $\gamma = 3/4$   
  &  80   &   3.0304e-03  &- &  2.8452e-03& - & 4.5423e-03& - \\    &  160   &   4.1974e-04 & 2.85& 3.3414e-04& 3.09&7.6124e-04&2.57\\
 &  320   &      5.9084e-05  & 2.82 &4.1481e-05 & 3.00&1.7839e-04&2.09\\
 &  640   &   7.5253e-06  & 2.97  &5.1517e-06& 3.00& 2.4096e-05 & 2.89 \\
 \hline
SI-R, $\gamma = 1-1/\sqrt{2}$                &  80   &  2.9584e-03  &- & 2.7743e-03& - & 4.4433e-03& - \\    &  160   &  4.0980e-04 & 2.85& 3.2398e-04& 3.09& 7.5081e-04&2.56\\
 &  320   &    5.7820e-05  & 2.83 &  4.0096e-05 & 3.01&1.7626e-04&2.09\\
 &  640   &   7.3666e-06  & 2.97  &4.9760e-06& 3.01&  2.3916e-05& 2.88 \\
 \hline
\end{tabular}
\end{center}
\caption{Test 2. The $L^2$, $L^1$, $L^{\infty}$ errors and orders of accuracy for equation (\ref{KdV_22}) with $\lambda=10$ and $T = \pi/4$.}
\label{table:KdV_10}
\end{table}%
\item[{Test 3.}] We conclude this section considering the following nonlinear biharmonic-type equation \cite{tan2022stability}
\begin{equation}\label{FourthO}
u_t  + ((u^2+2)u_{xx})_{xx} = f(x,t), \quad x \in (-\pi, \pi),
\end{equation}
with initial condition $u(x,0) =\sin(x)$, source term
\begin{equation}\label{f_sourse}
f(x,t) = e^{-3t}(e^{2t} - 6\cos^2(x) + 3 \sin^2(x))\sin(x),
\end{equation}
and exact solution $u(x, t) = e^{-t}\sin(x)$. We compute the solution to $T = 1$ with the time step $\Delta t = \Delta x$. 
The numerical errors and orders of accuracy are listed in Table \ref{table:1fourth}. We can see that the schemes are stable and achieve the expected orders of accuracy.
  
\begin{table}[htp]
\begin{center}
\begin{tabular}{c|c|c|c|c|c|c|c}
\hline
Scheme     & N     & $L^2$-error & order & $L^1$-error & order & $L^{\infty}$-error & order  \\
\hline
SI-R, $\gamma = 3/10$ &  40   &  2.4253e-03 & - &      2.3884e-03  &   - & 2.7038e-03 &    -      \\
                                       &  80   &     3.1623e-04  & 2.94 & 3.1654e-04& 2.91& 3.3670e-04 & 3.00\\
                                       &  160   &     4.1681e-05 & 2.92& 4.1714e-05& 2.92&4.3549e-05&2.95\\
                                      &  320   &       5.3676e-06  & 2.96 & 5.3466e-06 & 2.96 &5.7551e-06&2.92\\
                                      &  640   &       6.8320e-07  & 2.97  & 6.7657e-07& 2.98& 7.6768e-07& 2.91 \\
\hline
SI-R, $\gamma = 3/4$   &  40   &  2.1189e-02 & - &    1.8986e-02   &   - & 2.8878e-02 &      -    \\
                                       &  80   &   3.3998e-03  & 2.63 & 3.0674e-03& 2.62& 5.1576e-03 & 2.47\\
                                       &  160   &  4.6211e-04 & 2.87& 4.2009e-04& 2.86& {6.8619e-04}&2.91\\
                                      &  320   &      5.6802e-05  & 3.02 & 5.3937e-05 & 2.96 &8.0083e-05&3.09\\
                                      &  640   &    7.1809e-06  & 2.98  & 7.046e-06& 2.94& 9.0132e-06& 3.15 \\
                                      \hline
SI-R, $\gamma = 1-1/\sqrt{2}$   &  40   & 2.2745e-03 & - &        2.2475e-03   &   - & 2.5267e-03 &    -      \\
                                                    &  80   &   3.082e-04  & 2.88& 3.0777e-04& 2.86& 3.3221e-04 &2.92\\
                                                   &  160   &   4.1406e-05 &2.89& 4.1176e-05&2.90& 4.3276e-05&2.94\\
                                                   &  320   &      5.3855e-06  &2.94 &5.3251e-06 & 2.95&5.9901e-06& 2.85 \\
                                                   &  640   &  6.8859e-07  & 2.97  &  6.7672e-07& 2.98& 7.9635e-07& 2.91\\
                                                   \hline
\end{tabular}
\end{center}
\caption{Test 3. The $L^2$, $L^1$, $L^{\infty}$ errors and orders of accuracy for equation (\ref{FourthO}).}
\label{table:1fourth}
\end{table}%
\end{itemize}

\section{Partitioned semi-implicit Multistep methods.} \label{SIMM}
In \cite{albi2021high}, partitioned semi-implicit multistep methods were developed by applying distinct multistep schemes to the different components of system (\ref{ydotzdot}). 
Specifically, for the partitioned system (\ref{ydotzdot}), we consider schemes based on solving first component with an explicit linear multistep method and an implicit one for the second. For the sake of notation simplicity, we assume the system (\ref{ydotzdot}) is autonomous, such that the function $\mathcal{H}$ is independent of time $t$. Therefore, 
we obtain the semi-implicit linear multistep (SI-LM) solver as 
\begin{equation}\begin{aligned}\label{eq:semiGPS}
   u^{n+1} = - \sum_{j=0}^{\nuu-1} \tilde a_j u^{n-j} + \Delta t \sum_{j=0}^{\nuu-1} \tilde b_j {\mathcal H}\left(u^{n-j},v^{n-j}\right)\cr
   v^{n+1} = - \sum_{j=0}^{\nuu-1} a_j v^{n-j} +{\Delta t} \sum_{j=-1}^{\nuu-1} b_j  {\mathcal H}\left(u^{n-j},v^{n-j}\right),
\end{aligned}\end{equation} 
where $b_{-1} \neq 0$ and initially we assume $v^{n-j} = u^{n-j}$, $j=0,\ldots,\nuu-1$. Implicit linear multistep methods for which $b_j=0$, $j=0,\ldots,\nuu-1$ are referred to as backward differentiation formula (BDF). 
Another important class is represented by  Adams methods, for which $\tilde a_0=-1$, $a_0=-1$, $\tilde a_j=0$, $a_j=0$, $j=1,\ldots,\nuu-1$.  
Let us point out that even if the scheme doubles the number of unknown, the number of evaluations of the function ${\mathcal H}(u(t),v(t))$ is not doubled since both schemes use the same time levels. 

Now, based on equation (\ref{DiscS}), the function $\mathcal{H}$ has a linear dependence on its second argument, as follows:
\begin{align}
{\mathcal H}\left(u(t),v(t)\right)={\mathcal{F}}(u(t))+\mathcal{B}(u(t))v(t).
\label{eq:lin}
\end{align}
This particular choice allows for the resulting scheme to be solved without requiring an iterative solver.
 In fact, the second equation in \eqref{eq:semiGPS} can be rewritten as 
\begin{align*}
 v^{n+1} = 
 - \sum_{j=0}^{\nuu-1} a_j v^{n-j} +{\Delta t} \sum_{j=0}^{\nuu-1} b_j  {\mathcal H}\left(u^{n-j},v^{n-j}\right) 
 +{\Delta t} b_{-1}\left({\mathcal F}\left(u^{n+1}\right)+ {\mathcal B}(u^{n+1})v^{n+1}\right).
 \end{align*}
 or equivalently in explicit form
\begin{align}\label{Stab_a}
\begin{split}
 v^{n+1} =& \left(I-{b_{-1}\Delta t} {\mathcal B}(u^{n+1})\right)^{-1}\left(-\sum_{j=0}^{\nuu-1} a_j v^{n-j} + {\Delta t} \sum_{j=0}^{\nuu-1} b_j  {\mathcal H}\left(u^{n-j},v^{n-j}\right)\right.\\
 &+{\Delta t} b_{-1}\left({\mathcal F}\left(u^{n+1}\right)\right)\Bigg),
 \end{split}
 \end{align}
since $u^{n+1}$ is computed from the first equation in \eqref{eq:semiGPS}.

At variance with semi-implicit Runge-Kutta methods, 
multistep semi-implicit
scheme provides two distinct numerical solutions $u^{n+1}$ and
$v^{n+1}$ approximations of the true solution $u(t^{n+1})$ of problem
\eqref{ydotzdot}. Note, however, that both these solutions are used by the scheme to advance in time. {In order to define a unique solution}, it is natural to consider the scheme \eqref{eq:semiGPS} as a predictor--corrector multistep method for the non stiff component, {where an explicit scheme can be used} to predict $u^{n+1}$ which is then used by the implicit solver as a corrector for $v^{n+1}$.

\begin{rem}\label{Cc}
As a consequence of the above predictor-corrector formulation for the non stiff component, if the implicit solver has order $p$ it is typically enough to consider an explicit solver of order $p-1$, \cite{albi2021high}.
\end{rem}
    
However, the semi-implicit approach proposed in the paper \cite{albi2021high}, following the idea in remark (\ref{Cc}), does not ensure, in general,  the stability of the method. Specifically, this lack of stability arises from the explicit predictor step, which becomes problematic when the spatial discretization involves higher-order derivatives (e.g., a fourth derivative), thereby increasing the stiffness of the system.

For instance, we consider problem (\ref{DiscS}) with (\ref{eq:lin}) and set $\mathcal{F}(u(t)) = 0$. We choose the simplest explicit predictor with $p=1$ (namely, the explicit Euler scheme), 
and starting with $u^n = v^n, $we get for the numerical solution
\begin{align}
u^{n+1} &= u^n +   \frac{\Delta t}{ \Delta x^k}\mathcal{B}(u^n) v^{n} = \left(I + \frac{\Delta t}{ \Delta x^k}\mathcal{B}(u^n)\right)u^{n}.
\end{align}
\\
This equation involves an explicit evaluation of $u^{n+1}$, making the scheme conditionally stable, i.e.,
$$
\| u^{n+1}\| \le \| u^n\|,
$$ 
for a certain (semi-)norm, provided that a suitable time step restriction $\Delta t \le \Delta t_0$ is satisfied. For instance, in this case, we can choose $\Delta t_0 = C\Delta x^k/\left(\max_{u^n}(\mathcal{B}(u^n))\right)$.  
This implies that a stability constraint is imposed on the time step, which becomes increasingly restrictive as $k>2$, due to its dependence on $\Delta x^k$.\\  
However, this suggests that, in order to restore stability in the predictor step, and thereby making the scheme unconditionally stable, one can compute an $u^*$ using a semi-implicit approach of the form
\begin{align}\label{ImE}
u^{*} &= v^n +  \frac{\Delta t}{ \Delta x^k}\mathcal{B}(u^n) u^{*}. 
\end{align}
\\
This scheme corresponds to a backward Euler step and is unconditionally stable, in the sense that $\|u^{*}\| \le \|u^n\|$ for all $n$ and for any positive time-step $\Delta t > 0$. Indeed, the backward Euler method in \eqref{ImE} is unconditionally stable if the spectral radius $\rho(B)$ of the iteration matrix $B = \left(I -  \frac{\Delta t}{ \Delta x^k}\mathcal{B}(u^n) \right)^{-1}$ satisfies $\rho(B) < 1$. This condition holds whenever $\mathcal{B}(u^n)$ is nonsingular and negative definite (see, for instance, \cite{sebastiano2023high} for details). 
\\
Therefore, we will employ the SI approach even in the predictor step to restore the stability of the method, using a suitable SI predictor solver of order $p-1$ together with an SI corrector of order $p$. In the next section, we propose a {\em modified} SI-LM method to ensure both stability and accuracy.
\\
 \subsection{Semi-implicit predictor-corrector LM methods of order $p$}\label{Sec:SILM}
For the purpose of this section, assuming that the past values $u^{n-s}, ..., u^n$ are available, and  suppose that $v^{n-j} = u^{n-j}$, for $j = 0, ..., s-1$, we shall define a modified $s$-step SI-LM scheme. In practice, its implementation proceeds through the following steps:
\begin{align}\label{Algor}
\begin{split}
\textrm{predictor (P):} \\ 
   u^{*}  + \sum_{j=0}^{s-1} \bar a_j u^{n-j} &= \Delta t \sum_{j=0}^{s-1} \bar b_j \HH\left(u^{n-j-1},v^{n-j}\right) 
   + {\Delta t}  \bar{b}_{-1}  \HH\left( u^{n}, u^{*}\right),\\[2mm]
    &\textrm{(SI-LM predictor of order $p^*$)},\\[2mm]
\textrm{corrector (C):} \\
   v^{n+1} + \sum_{j=0}^{s-1} a_j v^{n-j} &= {\Delta t} \sum_{j=0}^{s-1} b_j  \HH\left(u^{n-j},v^{n-j}\right) 
   + {\Delta t}  b_{-1} \HH\left(u^{*},v^{n+1}\right),\\[2mm]
   &\textrm{(SI-LM corrector of order $p$)}
   \end{split}
\end{align}
Finally we update the solution $u^{n+1}= v^{n+1}$.\\
The coefficients  $\bar a_j^{(k)}$, $\bar b_j^{(k)}$, depend on the specific linear multistep scheme of order $p^*$ in the predictor step and $a_j$ and $b_j$ depend on the specific linear multistep scheme of order $p$ for the corrector one.
\\
Compared to implicit methods, in the semi-implicit predictor–corrector (SI-PC) approach (\ref{Algor}), a more practical procedure is usually to specify in advance how many iterations of the corrector are needed at each step to achieve the desired order of convergence. In this way, the scheme (\ref{Algor}) is both accurate and stable. We will apply the SI-PC pair (\ref{Algor}) in order to correct the order of convergence. 
Usually, the number of iterations of the corrector depends on the difference between orders $p^*$ and $p$.\\
In order to validate that our SI-PC scheme gets the corrector order, we start to introduce the local truncation error (LTE).
Following \cite{lambert1991numerical}, for a given linear multistep method
\begin{equation}\label{eq:lmm}
    v^{n+1} = - \sum_{j=0}^{\nuu-1} a_j v^{n-j} 
    + {\Delta t} \sum_{j=-1}^{\nuu-1} b_j  {\mathcal H}\left(t^{n-j},u^{n-j},v^{n-j}\right)
\end{equation}
we introduce the {\em linear difference operator} $\mathcal{L}$ associated with (\ref{eq:lmm}), defined by
\begin{equation}\label{LTEp}
\mathcal{L}(z(t);{\Delta t})=\sum_{j=0}^{k}\left[a_j z(t+j{\Delta t})-{\Delta t}b_jz'(t+j{\Delta t}))\right],
\end{equation}
where $z(t) \in C^1[a,b]$ is an arbitrary function. 
Hence for an order $p$ method, it holds:
$$\mathcal{L}(z(t); \Delta t) = C_{p+1}\Delta t^{p+1}z^{(p + 1)}(t) + \mathcal{O}(\Delta t^{p + 2}).$$
Let the predictor and corrector defined in  (\ref{Algor}) have associated linear difference 
operator $\mathcal{L}^*$ and $\mathcal{L}$, of order $p^*$, $p$ and error constant $C^*_{p^*+1}$ and $C_{p+1}$, respectively.
\\
Now we assume that all the past values are exact, that is, 
\begin{equation}\label{LA}
u^{n-j} = v^{n-j} = u(t_{n-j}),
\end{equation}
for $j =0,1,..,s$. Let us denote by ${u}^*_{n+1}$ and ${v}_{n+1}$ the values at time $t_{n+1}$ generated by the method when condition \eqref{LA} is in force.
 We assume that $u(t), v(t) \in \textrm{C}^{\hat{p}}[a,b]$,  where $\hat{p} = \max(p^*,p)$. 
Then one has for the LTE
\begin{align}\label{LTEs2}
\begin{split}
\mathcal{L}^{*}(u(t_n); \Delta t) &= C_{p^*+1}^*\Delta t^{p^*+1}u^{(p^* + 1)}(t_n) + \mathcal{O}(\Delta t^{p^* + 2}),\\
\mathcal{L}(v(t_n); \Delta t)   &= C_{p+1}\Delta t^{p+1}v^{(p + 1)}(t_n) + \mathcal{O}(\Delta t^{p + 2}).
\end{split}
\end{align}
From \eqref{LTEp}, we have for the predictor step
\[
u(t_{n+1}) + \sum_{j = 0}^{s-1}\bar a_ju(t_{n-j})  = \Delta t \sum_{j=0}^{s-1} \bar b_j \HH\left(u(t_{n-j-1}),v(t_{n-j})\right) 
   + {\Delta t}  \bar{b}_{-1}  \HH\left(u(t_n), u(t_{n+1}) \right) + \mathcal{L}^*(u(t_n);\Delta t).
\]
and, similarly, for the corrector one
\[
v(t_{n+1}) + \sum_{j = 0}^{s-1} a_j v(t_{n-j})  = \Delta t \sum_{j=0}^{s-1} b_j \HH\left(u(t_{n-j}),v(t_{n-j})\right) 
   + {\Delta t} {b}_{-1}  \HH\left( u(t_{n+1}), v(t_{n+1})\right) + \mathcal{L}(v(t_n);\Delta t).
\]
The values ${u}^*_{n+1}$ and ${v}_{n+1}$ given by the method satisfy
\[
u^*_{n+1} + \sum_{j = 0}^{s-1}\bar a_j u_{n-j}  = \Delta t \sum_{j=0}^{s-1} \bar b_j \HH\left(u_{n-j},{v}_{n-j}\right) 
   + {\Delta t}  \bar{b}_{-1}  \HH\left( u_{n},u^*_{n+1}\right),
\]
and
\begin{equation}\label{vns}
{v}_{n+1} + \sum_{j = 0}^{s-1} a_j {v}_{n-j}  = \Delta t \sum_{j=0}^{s-1} b_j \HH\left(u_{n-j},{v}_{n-j}\right) 
   + {\Delta t} {b}_{-1}  \HH\left(u^*_{n+1}, {v}_{n+1}\right).
\end{equation}
Subtracting and using the condition \eqref{LA} we obtain
\begin{equation}\label{uSol}
u(t_{n+1}) - {u}_{n+1} = \Delta t \bar{b}_{-1}\left(\mathcal{H}(u(t_n), u(t_{n+1})) - \mathcal{H}(u_n ,u^*_{n+1})\right) + \mathcal{L}^*(u(t_n);\Delta t),
\end{equation}
and 
\begin{equation}\label{vSol}
v(t_{n+1}) - {v}_{n+1} = \Delta t  b_{-1}\left(\mathcal{H}(u(t_{n+1}),v(t_{n+1})) - \mathcal{H}(u^*_{n+1},v_{n+1})\right) + \mathcal{L}(v(t_n);\Delta t).
\end{equation}
Now, from \eqref{LA} with $u(t_n) = u_n$, we have
$$
\mathcal{H}(u(t_n),u(t_{n+1})) - \mathcal{H}(u_{n},u^*_{n+1}) = \frac{\partial{\mathcal{H}}}{\partial u}[u(t_{n+1}) - u^*_{n+1}],
$$
where the derivative $\frac{\partial{\mathcal{H}}}{\partial u}$  denotes the Jacobian matrix of $\mathcal{H}$ with respect to $u$ evaluated at a point $\eta_{n+1}$ lying in the interior of the segment joining $u(t_{n+1})$ to $u^*_{n+1}$. 
Equation \eqref{uSol} now yields
\begin{equation}\label{uLTE}
[I - \Delta t\bar{b}_{-1} \partial{\mathcal{H}}/\partial u]
[u(t_{n+1}) - u^*_{n+1}] = \mathcal{L}^{*}(u(t_n);\Delta t). 
\end{equation}
Proceeding for the $v$-component, we have
\begin{equation}\label{vLTE}
[I - \Delta t\, b_{-1} \partial{\mathcal{H}}/\partial v]
[v(t_{n+1}) - v_{n+1}] = \Delta t \, b_{-1}\partial{\mathcal{H}}/\partial u \,[u(t_{n+1}) - u^*_{n+1}] + \mathcal{L}(v(t_n);\Delta t).
\end{equation}
In \eqref{uLTE} to a first-order approximation, that is, neglecting the $\mathcal{O}(\Delta t)$ term on the left-side, the LTE is simply the difference between the exact solution and the numerical one, subject to the condition \eqref{LA}, i.e.,
\begin{equation}\label{uLTE2}
u(t_{n+1}) - u^*_{n+1} = \mathcal{L}^{*}(u(t_n);\Delta t). 
\end{equation} 
Denoting by $\mu$ the number of correction steps, we first set $\mu = 1$. That is, we use the predictor value $u^*_{n+1}$ to compute $v_{n+1}$, and then evaluate the order of the SI–PC method after a single correction by analyzing the local truncation error LTE, under condition \eqref{LA}.
\\
Therefore, from \eqref{uLTE2}, for the $u$-component we have 
\begin{align}\label{LTEs}
u(t_{n+1}) - u^*_{n+1}  = C_{p^*+1}^*\Delta t^{p^*+1}u^{(p^* + 1)}(t_n) + \mathcal{O}(\Delta t^{p^* + 2}).
\end{align}
Whereas in order to evaluate the LTE for the $v$-component, we have to make the following considerations. 
\\
First consider the case $p^* \ge p$. On substituting \eqref{uLTE} into \eqref{vLTE} we get 
\begin{align}
v(t_{n+1}) - v_{n+1} = C_{p+1}\Delta t^{p+1}v^{(p + 1)}(t_n) + \mathcal{O}(\Delta t^{p + 2}).
\end{align}
Thus, if $p^*\ge p$ and $\mu = 1$, the LTE for the $v$-component is precisely the order of the corrector step.\\
Now we consider the case $p^* = p-1$ and denote by $v_{n+1}^{(\mu)}$ the numerical solution obtained from \eqref{vns} after the $\mu$-th correction step. For $\mu = 1$, on substituting again \eqref{uLTE} into \eqref{vLTE},  from \eqref{LTEs2} we obtain that
\begin{equation}\label{ConstantE}
v(t_{n+1}) - v_{n+1}^{(1)} =\left( b_{-1}\frac{\partial{\mathcal{H}}}{\partial u^*} C^*_p  u^{(p)}(t_n) + C_{p+1}v^{(p+1)}(t_n)\right)\Delta t^{p+1} + \mathcal{O}(\Delta t^{p+2}).
\end{equation}
Thus, after one correction $\mu = 1$, in the case $p^* = p-1$, the LTE of the SI-PC method is not identical with that of the corrector one even though it has the same order of accuracy as the corrector, since the leading error constant is different.\\
Now, setting $u^{*}_{n+1} = v_{n+1}^{(1)}$, with 
$u(t_{n+1}) - u^{*}_{n+1} = v(t_{n+1}) - v_{n+1}^{(1)}$,
we recompute with a second correction $\mu =2$, $v_{n+1}^{(2)}$ from \eqref{vns},
and we find for the LTE that 
$$
v(t_{n+1}) - v_{n+1}^{(2)} = \mathcal{L}(v(t_n); \Delta t) = C_{p+1}\Delta t^{p+1}v^{(p+1)}(t_n) + \mathcal{O}(\Delta t^{p+2}),
$$
namely, the SI-PC method has the same LTE of the corrector alone with the same leading error constant.\\
Finally, we consider the case $p^* = p-2$ and $\mu = 1$. On substituting \eqref{uLTE} into \eqref{vLTE}, and from \eqref{LTEs2}, we get that
\begin{equation}\label{ConstantE_2}
v(t_{n+1}) - v_{n+1}^{(1)} = b_{-1}\frac{\partial{\mathcal{H}}}{\partial u^*} C^*_{p-1}  u^{(p-1)}(t_n) \Delta t^{p} + \mathcal{O}(\Delta t^{p+1}).
\end{equation}
Thus, if $\mu = 1$ the order of the SI-PC method is only $p-1$. \\
Now, proceeding as before, we recompute $v_{n+1}^{(2)}$ from \eqref{vns} 
and on substituting \eqref{ConstantE_2} into \eqref{vLTE} we obtain an order $p+1$ for the LTE of the SI-PC method, which is identical to that of the corrector. Nevertheless, also in this case, the corresponding leading error constant is not the same. Continuing this process iteratively, for $\mu \ge 3$, and using again \eqref{vLTE} with $v_{n+1}^{(\mu)}$ 
the SI-PC method achieves the same LTE of the corrector alone. 
\\
We therefore state the following proposition.
\begin{prop}\label{Pro_1}
  Consider the differential equation (\ref{ydotzdot}), with $\mathcal{H}(t,u,v)$ continuously differentiable and let $v(t) = u(t)$ be its solution. Assume condition \eqref{LA} and consider a SI-PC method \eqref{Algor}, where the predictor method (P) has order $p^*$ and the corrector one has order $p$.  Then we have the following cases:
  \begin{itemize}
  \item if $p^*\ge p$, for all $\mu \le 1$ the LTE for the numerical solution $v_{n+1}$ is precisely that of the corrector (C) alone. 
    \item if $p^* = p-1$, with $\mu = 1$  the LTE of the SI-PC method \eqref{Algor} is not identical with that of the corrector even though it has the same order of accuracy as the corrector. With $\mu \ge 2$, the LTE for the numerical solution $v_{n+1}$ is precisely that of the corrector (C) alone.
    \item if $p^* < p$, with $\mu = p - 
    p^*$,  the SI-PC method and corrector have the same order but different error constants. With $\mu > p-p^*$ the order is that of the corrector step.
    \end{itemize}
 \end{prop} 
In order to derive our semi-implicit predictor-corrector method, 
we choose as a predictor step the SI-Euler method with $p^* = 1$ and correct the numerical solution $v_{n+1}^{(\mu)}$ in \eqref{vns} by a fixed number  $\mu$ of suitable iterations %
in order to reach the correct order of convergence. For instance,
let us choose a SI-BDF method of order $p$ for the correction steps in \eqref{vns}. We call this type of predictor-corrector schemes SI-PC$^{\mu}$ BDF$p$.
\\ 
\begin{example}\label{example:si-pc-bdf}
For example if we use a SI-BDF3 as corrector,
in order to achieve the correct order of convergence, i.e. $p = 3$,  we require at least $\mu = 2$. 
Note that with $\mu=3$ we can improve the accuracy because the third iteration also allows one to recover the error constant of the corrector method, while $\mu>3$ 
does not provide further improvement.\\
Below we provide the algorithm with SI-Euler method as predictor and SI-BDF3 as corrector, in order to achieve the third order with $\mu = 2$ (or $\mu = 3$).\\
\noindent\textbf{Algorithm: SI-PC$^{\,\mu}$ BDF3}
\begin{align} \label{algor}
\begin{split}
&\textbf{Predictor (P):} \\
&\quad u^{*,(0)} = u^n + \Delta t\, \mathcal{H}(u^n, u^{*,(0)}),
\qquad \text{(SI-Euler, order } p = 1\text{)} \\[2mm]
&\textbf{Corrector (C):} \\[-1mm]
&\quad \text{for } \nu = 0, 1, ..., \mu-1 \text{ do} \\[1mm]
&\qquad u^{*} = u^{*,(\nu)}, \\[1mm]
&\qquad u^{*,(\nu + 1)} = 
\frac{18}{11} u^n - \frac{9}{11} u^{n-1} + \frac{2}{11} u^{n-2}
+ \frac{6}{11}\, \Delta t\, \mathcal{H}(u^{*}, u^{*,(\nu + 1)}), \\[1mm]
&\qquad \text{(SI-BDF3, order } p = 3\text{)} \\[1mm]
&\quad \text{end for} \\[2mm]
&\textbf{Set:} \quad v^{n+1} = u^{*,(\nu + 1)}.
\end{split}
\end{align}    
\end{example}

In the numerical tests section, we have chosen BDF methods as our linear multistep schemes because of their good stability properties. \\
In \eqref{algor}, to perform the predictor step and the $\mu$ correction steps, a certain number of linear systems must be solved.  
Specifically, one linear system is required for the predictor (the SI-Euler method) and $\mu$ systems for the corrections, resulting in a total of $\mu + 1$ linear systems per time step.\\
For instance, to achieve third-order accuracy in \eqref{algor}, three linear systems must be solved. However, by solving a fourth linear system, one can get the same error constant of the corrector method.
Similarly, if the corrector is replaced by a SI-BDF4 scheme, to obtain fourth-order accuracy at least four linear systems are required. In general, by means of this strategy, a BDF of order $p$ requires solving $p $ or $p + 1$ linear systems per time step.
This behavior will be illustrated in the numerical tests presented 
in Section \ref{NumEx}.\\
In any case, once we choose an order $p$ corrector method, we could choose any suitable order $p-1$ predictor method, such as for example
a {\em chain} of semi-implicit linear multistep methods of order $1$ to $p-1$. Of course, this choice can affect the computational costs of the overall method.
\subsection{Computational Aspects.} 
In order to show the effective benefit of these methods, we remark again that in general these semi-implicit predictor-corrector schemes provide 
a way to efficiently compute an accurate solution 
of problem (\ref{G-systemBIS}) by means of a fixed number of linear systems solve, even when the problem presents stiffness
due to high order spatial derivatives.\\
Indeed, a semi-implicit predictor-corrector BDF scheme of order $p$ requires solving $p$ or $p + 1$ linear systems per time step, 
and does not suffer from instability as in the case when using an explicit predictor.

The semi-implicit predictor-corrector BDFs have computational costs per step
comparable to SI-RK methods \cite{sebastiano2023high} that have order $p\leq3$, 
but provide a simple and effective way to obtain higher order schemes,
since they are based on regular linear multistep methods. 
Indeed, semi-implicit RK schemes are based on existing additive IMEX methods \cite{boscarino2016high}, 
so, for methods of order $p\ge 3$ are needed at least $p+1$ stages and hence $p+1$ linear systems to solve.
Furthermore, only a few high-order (order $p > 3$) additive IMEX-RK schemes exist in the literature (e.g.,
Carpenter et al., \cite{kennedy2019higher}), and they typically require 6 to 8 stages to reach fourth- or fifth-order
accuracy. Finally, existing IMEX RK schemes have low stage order, and hence they can suffer from order reduction.

In any case, since the proposed methods are multistep, a suitable {\em starting procedure} is needed \cite{lambert1991numerical, Calab2023, Izzo2025}.
As usual in the context of multistep linear methods, 
the starting values can be approximated by means of a
SI-PC$^\mu$ BDF method of order $p-1$ applied with a reduced
stepsize $\Delta t/m$, with an integer, suitably chosen $m$.
In the numerical tests here reported a maximum value of $m=4$ has been used for Test 1 and 2, while a maximum value of $m=16$ has been used for Test 3 and 4.

\subsection{Numerical Results}\label{NumEx}
In this section, we will numerically verify the orders of accuracy and performance of SI-PC multistep schemes \eqref{Algor} for the solution of time dependent PDEs with high order derivatives, as equations: (\ref{First1}), (\ref{Conv-Diff}),  (\ref{eq:Disper2_bis}), (\ref{eq:Fourth}) and (\ref{eq:Fourth2}). For simplicity we used in all of our examples periodic boundary conditions.

In each of the numerical tests presented below, we consider second, third and fourth order modified SI-PC multistep schemes described in Section \ref{Sec:SILM} and in Example \ref{example:si-pc-bdf}: SI-PC$^2$ BDF2, SI-PC$^3$ BDF3 and SI-PC$^4$ BDF4 for the time discretization. We will provide the time step $\Delta t$ that will produce stable solutions across the entire solution domain. To achieve order $p$ with SI-PC$^\mu$ BDF$p$, we have experimentally observed that solving only $p$ linear systems per
time step, may often be sufficient. In some cases, depending on the problem, one additional correction step may be needed to reach
the desired accuracy. For the sake of brevity and presentation, we have not reported here the results obtained by SI-PC$^\mu$ BDF$p$ with $\mu=p-1$. 

For the space discretization we use high-order finite difference schemes 
\cite{shu2006essentially}, in particular, in order to reach the global second,  third and four order accuracy we consider  the standard WENO(3,2) and WENO(5,3) schemes \cite{shu2006essentially} for the discretization of the convection term presented in the following tests.

\begin{itemize}
\item[Test 1.] First we consider the diffusion equation
\begin{equation}\label{Eq:1Diff}
u_t = (a(u)u_x)_x + f(x,t), \quad  x\in (-\pi, \pi),
\end{equation}
Here we consider the case $a(u) = u^2 +1$ with initial condition $u(x,0) = \sin(x)$ and the source term $f(x,t)$ is chosen such that the exact solution is $u(x,t) = \sin(x-t)$. For this example, we take: the time step $\Delta t = \Delta x$, and $\Delta x = 2\pi/N$, with  periodic boundary conditions and final time $T = 10$. 
In Table \ref{table:1SILM} we display the numerical results of 
SI-PC$^\mu$ BDF$p$ schemes, with $\mu = p$ 
and $p=2,3,4$, coupled with the finite difference type spatial discretization. 
In the table we list the $L^2$, \, $L^1$ and $L^{\infty}$ errors and the orders of accuracy of each scheme. We can see that the schemes are stable and achieve the expected order of accuracy.\\

\begin{table}[!htbp] 
\begin{center} 
\begin{tabular}{c|c|c|c|c|c|c|c} 
\hline 
Scheme     & N     & $L^2$-error & order & $L^1$-error & order& $L^\infty$-error & order \\ [1mm] 
\hline 
 SI-PC$^2$ BDF2 & 40 & 4.9945e-03 &  -    & 4.8052e-03 &  -    & 5.6425e-03 &  -   \\ 
   & 80 & 1.2655e-03 & 1.98  & 1.2155e-03 & 1.98  & 1.4290e-03 & 1.98 \\ 
   & 160 & 3.2178e-04 & 1.98  & 3.0859e-04 & 1.98  & 3.6369e-04 & 1.97 \\ 
   & 320 & 8.1720e-05 & 1.98  & 7.8267e-05 & 1.98  & 9.2704e-05 & 1.97 \\ 
\hline 
 SI-PC$^3$ BDF3 & 40 & 5.1979e-04 &  -    & 5.3682e-04 &  -    & 4.7807e-04 &  -   \\ 
   & 80 & 6.3538e-05 & 3.03  & 6.4531e-05 & 3.06  & 6.2774e-05 & 2.93 \\ 
   & 160 & 7.9072e-06 & 3.01  & 7.9777e-06 & 3.02  & 8.1238e-06 & 2.95 \\ 
   & 320 & 9.8617e-07 & 3.00  & 9.9384e-07 & 3.00  & 1.0236e-06 & 2.99 \\ 
\hline 
 SI-PC$^4$ BDF4 & 40 & 1.0994e-04 &  -    & 1.0398e-04 &  -    & 1.2210e-04 &  -   \\ 
   & 80 & 6.9958e-06 & 3.97  & 6.5808e-06 & 3.98  & 7.8412e-06 & 3.96 \\ 
   & 160 & 4.3934e-07 & 3.99  & 4.1397e-07 & 3.99  & 5.0557e-07 & 3.96 \\ 
   & 320 & 2.7497e-08 & 4.00  & 2.5921e-08 & 4.00  & 3.1928e-08 & 3.99 \\ 
\hline 
\end{tabular} 
\caption{The $L^2$, $L^1$, $L^{\infty}$ errors and orders of accuracy for the second order diffusion equations (Test 1), for SI-PC$^\mu$ BDF$p$ methods with $\mu=p$ iterations and $p = 2,3,4$. \label{table:1SILM}} 
\end{center} 
\end{table} 

Next, we consider the same numerical tests introduced in section \eqref{sect:NR} and we show that all the schemes are stable and achieve the correct order of accuracy.

\item[Test 2.] The convection-diffusion equation \eqref{Eq:2Diff} with the same initial condition, source term and exact solution.  
Final time $T = 4$ and time step $\Delta t =  C\Delta x$ with $C = 4$. The numerical results are listed in Table \ref{table:alpha}.

\begin{table}[!htbp] 
\begin{center} 
\begin{tabular}{c|c|c|c|c|c|c|c} 
\hline 
Scheme     & N     & $L^2$-error & order & $L^1$-error & order& $L^\infty$-error & order \\ [1mm] 
\hline 
  SI-PC$^2$ BDF2 & 40 & 4.7704e-02 &  -    & 4.7454e-02 &  -    & 5.2723e-02 &  -   \\ 
   & 80 & 1.3009e-02 & 1.87  & 1.2740e-02 & 1.90  & 1.5084e-02 & 1.81 \\ 
   & 160 & 3.3349e-03 & 1.96  & 3.2462e-03 & 1.97  & 3.9278e-03 & 1.94 \\ 
   & 320 & 8.4052e-04 & 1.99  & 8.1573e-04 & 1.99  & 9.9778e-04 & 1.98 \\ 
\hline 
 SI-PC$^3$ BDF3 & 40 & 2.2582e-02 &  -    & 2.2040e-02 &  -    & 2.6428e-02 &  -   \\ 
   & 80 & 2.7192e-03 & 3.05  & 2.6847e-03 & 3.04  & 3.1086e-03 & 3.09 \\ 
   & 160 & 3.3292e-04 & 3.03  & 3.2998e-04 & 3.02  & 3.8601e-04 & 3.01 \\ 
   & 320 & 4.1807e-05 & 2.99  & 4.1239e-05 & 3.00  & 5.0343e-05 & 2.94 \\ 
\hline 
 SI-PC$^4$ BDF4 & 40 & 1.0584e-02 &  -    & 1.0858e-02 &  -    & 1.0397e-02 &  -   \\ 
   & 80 & 7.2030e-04 & 3.88  & 7.1663e-04 & 3.92  & 7.9986e-04 & 3.70 \\ 
   & 160 & 4.6143e-05 & 3.96  & 4.5104e-05 & 3.99  & 5.2604e-05 & 3.93 \\ 
   & 320 & 2.8343e-06 & 4.03  & 2.7773e-06 & 4.02  & 3.1524e-06 & 4.06 \\ 
\hline 
\end{tabular} 
\caption{The $L^2$, $L^1$ and $L^{\infty}$ errors and orders of accuracy for equation (\ref{Eq:2Diff})  (Test 2) 
for SI-PC$^\mu$ BDF$p$ methods with $\mu=p$ iterations and $p = 2,3,4$.
\label{table:alpha}}
\end{center} 
\end{table}

\item[Test 3.] 
Test problem (\ref{KdV_22}), i.e., the general KdV equation \cite{tan2022stability}, with the same initial condition and exact solution. 
Final time $T = \pi$, with ${\lambda} = 0.1$ and 
we set the time step as \eqref{Delta_t} with $CFL=0.4$.
In Table \ref{table:alpha2} are shown the results.

\begin{table}[!htbp] 
\begin{center} 
\begin{tabular}{c|c|c|c|c|c|c|c} 
\hline 
Scheme     & N     & $L^2$-error & order & $L^1$-error & order& $L^\infty$-error & order \\ [1mm] 
\hline 
 SI-PC$^2$ BDF2 & 40 & 6.1589e-03 &  -    & 5.5845e-03 &  -    & 1.2083e-02 &  -   \\ 
   & 80 & 8.3338e-04 & 2.89  & 7.3968e-04 & 2.92  & 1.3829e-03 & 3.13 \\ 
   & 160 & 1.0540e-04 & 2.98  & 9.1965e-05 & 3.01  & 2.6687e-04 & 2.37 \\ 
   & 320 & 1.3140e-05 & 3.00  & 1.1300e-05 & 3.02  & 4.4219e-05 & 2.59 \\ 
\hline 
 SI-PC$^3$ BDF3 & 40 & 6.1748e-03 &  -    & 5.5947e-03 &  -    & 1.2194e-02 &  -   \\ 
   & 80 & 8.3732e-04 & 2.88  & 7.4289e-04 & 2.91  & 1.3927e-03 & 3.13 \\ 
   & 160 & 1.0633e-04 & 2.98  & 9.2869e-05 & 3.00  & 2.6703e-04 & 2.38 \\ 
   & 320 & 1.3406e-05 & 2.99  & 1.1534e-05 & 3.01  & 4.4895e-05 & 2.57 \\ 
\hline 
 SI-PC$^4$ BDF4 & 40 & 5.3594e-04 &  -    & 4.2829e-04 &  -    & 1.2219e-03 &  -   \\ 
   & 80 & 3.2804e-05 & 4.03  & 2.5719e-05 & 4.06  & 1.0037e-04 & 3.61 \\ 
   & 160 & 2.0451e-06 & 4.00  & 1.6402e-06 & 3.97  & 6.9656e-06 & 3.85 \\ 
   & 320 & 1.3039e-07 & 3.97  & 1.0193e-07 & 4.01  & 5.0042e-07 & 3.80 \\ 
\hline 
\end{tabular} 
\caption{The $L^2$, $L^1$ and $L^{\infty}$ errors and orders of accuracy for equation (\ref{KdV_22})  (Test 3) 
for SI-PC$^\mu$ BDF$p$ methods with $\mu=p$ iterations and $p = 2,3,4$.
\label{table:alpha2}} 
\end{center}
\end{table}

\item[Test 4.] The nonlinear biharmonic-type equation \eqref{FourthO}, 
with initial condition $u(x,0) =\sin(x)$, source term
$f(x,t)$ in \eqref{f_sourse} 
and exact solution $u(x, t) = e^{-t}\sin(x)$.
Final time $T = 1$ and time step $\Delta t = \Delta x$. 
The numerical errors and orders of accuracy are listed in {Tables  \ref{Tab:4}.

\begin{table}[!htbp] 
\begin{center} 
\begin{tabular}{c|c|c|c|c|c|c|c} 
\hline 
Scheme     & N     & $L^2$-error & order & $L^1$-error & order& $L^\infty$-error & order \\ [1mm] 
\hline 
 SI-PC$^2$ BDF2 & 40 & 4.1028e-03 &  -    & 4.1125e-03 &  -    & 4.0727e-03 &  -   \\ 
   & 80 & 1.0353e-03 & 1.99  & 1.0377e-03 & 1.99  & 1.0282e-03 & 1.99 \\ 
   & 160 & 2.6795e-04 & 1.95  & 2.6849e-04 & 1.95  & 2.6632e-04 & 1.95 \\ 
   & 320 & 6.7449e-05 & 1.99  & 6.7583e-05 & 1.99  & 6.7045e-05 & 1.99 \\ 
\hline 
 SI-PC$^3$ BDF3 & 40 & 4.5862e-04 &  -    & 4.5917e-04 &  -    & 4.5693e-04 &  -   \\ 
   & 80 & 5.9486e-05 & 2.95  & 5.9583e-05 & 2.95  & 5.9193e-05 & 2.95 \\ 
   & 160 & 7.6394e-06 & 2.96  & 7.6525e-06 & 2.96  & 7.6002e-06 & 2.96 \\ 
   & 320 & 9.4389e-07 & 3.02  & 9.4559e-07 & 3.02  & 9.3887e-07 & 3.02 \\ 
\hline 
 SI-PC$^4$ BDF4 & 40 & 3.5466e-05 &  -    & 3.5876e-05 &  -    & 3.4130e-05 &  -   \\ 
   & 80 & 2.8023e-06 & 3.66  & 2.8302e-06 & 3.66  & 2.7137e-06 & 3.65 \\ 
   & 160 & 1.8669e-07 & 3.91  & 1.8835e-07 & 3.91  & 1.8149e-07 & 3.90 \\ 
   & 320 & 1.2020e-08 & 3.96  & 1.2124e-08 & 3.96  & 1.1703e-08 & 3.95 \\ 
\hline 
\end{tabular} 
\caption{The $L^2$, $L^1$, $L^{\infty}$ errors and orders of accuracy for 
fourth order diffusion equation  (Test 4),
for SI-PC$^\mu$ BDF$p$ methods with $\mu=p$ iterations and $p = 2,3,4$.
}
\label{Tab:4}
\end{center} 
\end{table} 

}
\end{itemize}

\section{Concluding remarks}
In this paper, following the ideas developed in \cite{sebastiano2023high} we have developed a SI strategy for Rosenbrock type and multistep methods based on an implicit-explicit setting for solving high order dissipative, dispersive and special biharmonic-type equations in one dimension. The SI schemes so designed for high order time-dependent PDEs 
do not require nonlinear iterative solvers and avoid the severe time-step restrictions
typically encountered when using explicit methods. Numerical experiments show that the schemes are stable and achieve the expected orders of accuracy for large time step.
 
 Since the main issue of this paper is to design an appropriate high-order time discretization for time-dependent PDEs, the numerical experiments are only performed for one-dimensional equations. For the space discretization, we considered only classical finite difference spatial discretization because  of its simplicity in design and coding and it is straightforward to extend to higher-dimensional equations. However, other types of space discretization can be used.

 \section*{Acknowledgements}
Sebastiano Boscarino is supported for this work by 1) the Spoke 1 "FutureHPC $\&$ BigData" of the Italian Research Center on High-Performance Computing, Big Data and Quantum Computing (ICSC) funded by MUR Missione 4 Componente 2 Investimento 1.4: Potenziamento strutture di ricerca e creazione di "campioni nazionali di R$\&$S (M4C2-19 )", (CN00000013);
by 2) the Italian Ministry of Instruction, University and Research (MIUR) to support this research with funds coming from PRIN Project 2022  (2022KA3JBA), entitled "Advanced numerical methods for time dependent parametric partial differential equations and applications"; 3) from Italian Ministerial grant PRIN 2022 PNRR "FIN4GEO: Forward and Inverse Numerical Modeling of hydrothermal systems in volcanic regions with application to geothermal energy exploitation.", (No. P2022BNB97). 
\\
Giuseppe Izzo is supported by the Italian MUR under the PRIN 2022 project No. 2022N3ZNAX.
\\
Sebastiano Boscarino and Giuseppe Izzo are members of the INdAM Research group GNCS.


\vspace{5mm}
\noindent {\bf Conflict of interest}\\
\noindent The authors declare that they have no conflict of interest.

\vspace{5mm}
\noindent {\bf Data availability} \\
Data sharing not applicable to this article as no datasets were generated or analysed during the current study.

%

\bibliographystyle{plain}
\bibliography{refs}

\end{document}